# DISCRETE LEAST-SQUARES FINITE ELEMENT METHODS

DEDICATED TO J. T. ODEN ON THE OCCASION OF HIS $80^{\text{TH}}$ BIRTHDAY.

BRENDAN KEITH, SOCRATIS PETRIDES, FEDERICO FUENTES, AND LESZEK DEMKOWICZ

ABSTRACT. A finite element methodology for large classes of variational boundary value problems is defined which involves discretizing two linear operators: (1) the differential operator defining the spatial boundary value problem; and (2) a Riesz map on the test space. The resulting linear system is overdetermined. Two different approaches for solving the system are suggested (although others are discussed): (1) solving the associated normal equation with linear solvers for symmetric positive-definite systems (e.g. Cholesky factorization); and (2) solving the overdetermined system with orthogonalization algorithms (e.g. QR factorization). The finite element assembly algorithm for each of these approaches is described in detail. The normal equation approach is usually faster for direct solvers and requires less storage. The second approach reduces the condition number of the system by a power of two and is less sensitive to round-off error. The rectangular stiffness matrix of second approach is demonstrated to have condition number $\mathcal{O}(h^{-1})$ for a variety of formulations of Poisson's equation. The stiffness matrix from the normal equation approach is demonstrated to be related to the monolithic stiffness matrices of least-squares finite element methods and it is proved that the two are identical in some cases. An example with Poisson's equation indicates that the solutions of these two different linear systems can be nearly indistinguishable (if round-off error is not an issue) and rapidly converge to each other. The orthogonalization approach is suggested to be beneficial for problems which induce poorly conditioned linear systems. Experiments with Poisson's equation in single-precision arithmetic as well as the linear acoustics problem near resonance in double-precision arithmetic verify this conclusion. The methodology described here was developed as an outgrowth of the discontinuous Petrov-Galerkin (DPG) methodology of Demkowicz and Gopalakrishnan [29]. The strength of DPG is highlighted throughout, however, the majority of theory presented is general. Extensions to constrained minimization principles are also considered throughout but are not analyzed in experiments.

## 1. INTRODUCTION

Least-squares finite element methods have been demonstrated to be an auspicious class of methods for a wide variety of boundary value problems of engineering interest. These methods are attractive for many challenging problems because of their general simplicity of implementation, their adherence to simple minimization principles and, therefore, numerical stability, and their built-in *a posteriori* error estimator. As we demonstrate in Section 5.1, the monolithic coefficient matrix in a least-squares finite element method can be identified with the coefficient matrix in the normal equation corresponding to a feasible linear least-squares problem. For a thorough study of least-squares methods and the most significant references, we refer to Bochev and Gunzburger [11].

In this article, we present *discrete least-squares* (DLS) finite element methods. As the name suggests, they can be related to least-squares finite element methods, however, the principles governing DLS methods are much more general. Least-squares finite element methods always induce monolithic, Hermitian positive-definite stiffness matrices. On the other hand, DLS methods can generate an overdetermined system of equations, as well as a related normal equation with a Hermitian positive-definite stiffness matrix.

All DLS methods involve a minimum residual principle and a discretization procedure where two operators are discretized: (1) The differential operator defining the spatial boundary value problem; and (2) a test space Riesz map. As we show in Section 3, this induces a discrete minimum residual principle for the solution vector $\mathbf{u}$ which can be expressed as a least-squares problem of the form

$$\min_{\mathbf{u}\in\mathbb{F}^N} \|\mathsf{B}\mathbf{u}-\mathsf{l}\|_{\mathsf{G}^{-1}} \qquad \Longleftrightarrow \qquad \mathbf{u} = (\mathsf{B}^*\mathsf{G}^{-1}\mathsf{B})^{-1}\mathsf{B}^*\mathsf{G}^{-1}\mathsf{l}\,,$$

*Date*: May 4, 2017.

where the field $\mathbb{F}$ is either $\mathbb{R}$ or $\mathbb{C}$, $\mathsf{G} \in \mathbb{F}^{M\times M}$ is Hermitian positive (semi-)definite and $\mathsf{B} \in \mathbb{F}^{M\times N}$ is rectangular, $M > N$. Another useful perspective is that, by duality, all DLS methods can be identified with the special equilibrium system [68]:

$$\begin{bmatrix} \mathsf{G} & \mathsf{B} \\ \mathsf{B}^* & 0 \end{bmatrix} \begin{bmatrix} \mathbf{p} \\ \mathbf{u} \end{bmatrix} = \begin{bmatrix} \mathsf{l} \\ 0 \end{bmatrix} . \tag{1.1}$$

Here, $\mathbf{u}$ is the subvector of primary interest, $\mathsf{G}$ can often be chosen by the numerical analyst, and there is always the orthogonality condition in the dual variable: $\mathsf{B}^*\mathbf{p} = 0$. This leads to many special features of DLS not common in most finite element saddle-point problems.

Our principal examples will pertain to a specific class of DLS methods: discontinuous Petrov-Galerkin (DPG) methods [29]. In this setting, the matrix $\mathsf{G}$ can be made block-diagonal, no matter the order of the interpolation spaces used in discretization. It is this feature that makes DPG methods the most practical and general DLS setting we are presently aware of. However, because the discrete saddle-point problem (1.1) can appear in several non-DPG finite element settings [73, 24, 54, 20, 14, 55], it is appropriate to consider this entire class of methods from a comprehensive perspective.

In some DLS methods, both subvectors in (1.1) are be jointly solved for because $\mathsf{G}$ cannot be efficiently inverted [24, 20]. Because $\mathsf{G}$ *can* be efficiently inverted with DPG, in that setting, much smaller problems, posed solely in the primal variable $\mathbf{u}$, can be solved directly. This has been performed for many problems of engineering interest [65, 21, 67, 33, 37, 49, 62, 30, 40, 38, 35].

In the DPG literature, the minimum residual principle is usually introduced in an idealized scenario where the test space Riesz map is known exactly—that is, it has not been discretized—and from this preternatural position, many notable properties of the method are inferred. Only later, when presenting a computer implementation of the method, a tunable discretization of the test space is proposed with the assumption that all properties of the idealized method—such as stability—will be inherited once the test space discretization is rich enough. In some scenarios, these assumptions can be proven [45, 18, 19, 35, 40, 39] or analyzed [56], but for difficult problems, they can only be speculated. We will not devote our attention in this article to such implications of different discretizations, but rather deliberate the structure of these equations mainly from a numerical linear algebra perspective. Finally, when the Riesz map in a DPG method is discretized, this is referred to as the "practical" DPG method by the majority of the community. It is only in this sense that DPG can be considered a DLS method.

One advantage of the perspective we have taken is clear in Section 3.1, where DPG is identified with an overdetermined system of equations. It is probably because the computer implementation of the practical DPG method has only begun to be investigated in its own right for efficiency and accuracy in large-scale computations [64, 66, 4], that this connection with an overdetermined discrete least-squares problem was overlooked, in the past. Furthermore, by exposing the connection with DLS, this article demonstrates that the established manner of assembling a practical DPG problem can affect the round-off error in the solution in a substantial way.

*Critical assumptions.* An additional benefit of the generality of this work is its subtle suggestions at efficient possibilities for DLS outside the DPG setting. Indeed, even in a DPG method, the "ideal" Gram matrix $\mathsf{G}_{\text{ideal}}$, coming from an idealized norm, is usually substituted with a more easily computable Gram matrix $\mathsf{G}$ coming from a similar norm. If the norms are equivalent, then only the *converged* solution of the dual variable is affected by this modification of the test norm, and not the *converged* solution of the primal variable. Because the exact $\mathsf{G}_{\text{ideal}}^{-1}$ and the exact $\mathsf{G}^{-1}$ could vary greatly—indeed, one could be dense and the other sparse—it is possible that satisfactory estimates of $\mathsf{G}_{\text{ideal}}^{-1}$ may be obtainable without using a block diagonal $\mathsf{G}$.

Therefore, we state that the three principal scenarios where DLS is most appropriate are:

(1) When test functions are constructed so that they are orthogonal under the chosen test space norm. In this situation, the Gram matrix $\mathsf{G}$ would be diagonal and its inverse would be trivial to compute.
(2) When the test space is broken—as is the case with DPG methods—the Gram matrix $\mathsf{G}$ can be directly inverted, because the computation would be local and feasible.
(3) When the Gram matrix $\mathsf{G}$ has a sparse inverse [69], or it is only approximately inverted but its approximate inverse is sparse.

*Specific remarks for the DPG community.* As far as we, the authors, are aware, in all published papers on the DPG method/methodology, the so-called "practical" DPG method for ultraweak formulations (with broken test spaces) is identical to solving the normal equation of the DLS system presented within. There appear to be no published manuscripts which derive or consider the possible advantages inherent in solving the overdetermined system of equations.

We note that in specific instances where the test space is $L^2$ in at least one component—and therefore, for most alternative variational formulations [19, 48]—the implementation of DPG permits some ambiguity. For instance, as developed in [48] and contrary to the DLS approach, one can choose not to discretize the $L^2$-Riesz map and so construct a least-squares (strong variational formulations) or least-squares hybrid (all other variational formulations) discretization.

The authors believe that this ultimately leads to some confusion. Different discretizations will involve different analysis, different advantages and limitations, and ultimately different results, however slight. We argue to adopt a convention to distinguish between these approaches. That is, the discretization advocated for here should aptly be called the *discrete* least-squares discretization, and the least-squares and hybrid discretizations developed in [48] should be identified as such.

We further advocate that the "ideal" DPG methodology remain as the consummate—however, often impractical—method which derives from the union of a minimum residual principle with broken test spaces and a Riesz map which is manipulated *analytically*. With this distinction, the ideal DPG methodology applied to a strong formulation is identical to traditional least-squares formulations while in the DLS setting, it is a different and unique numerical method.

*Other pertinent remarks, definitions, and notation.* Since we are synthesizing ideas from different fields, we require a certain amnesty on behalf of the reader throughout this article. For example, let $\mathsf{M} \in \mathbb{F}^{M\times N}$ be full rank matrix, $M > N$, and $\mathsf{b} \in \mathbb{F}^M$. In the least-squares context, the Gram (or Gramian) matrix is usually defined to be the matrix $\mathsf{M}^*\mathsf{M}$ in the normal equation

$$\mathsf{M}^*\mathsf{M}\mathsf{x} = \mathsf{M}^*\mathsf{b}\,.$$

In this article, we will require significantly different notion of Gram matrix. Accordingly, because our perspective is towards finite element methods, we will generally refer to matrices of the form $\mathsf{M}^*\mathsf{M}$ as *stiffness matrices*. Analogously, rectangular matrices like $\mathsf{M}$ will be referred to as an *enriched stiffness matrices*.

For most methods we have in mind, the trial functions $\mathfrak{u} = (\mathfrak{u}_1, \ldots, \mathfrak{u}_n)$ and test functions $\mathfrak{v} = (\mathfrak{v}_1, \ldots, \mathfrak{v}_m)$, have many components. The boundary value problems will be posed on simply connected Lipschitz domains $\Omega \subset \mathbb{R}^d$, $d > 0$, with sufficiently smooth boundaries $\partial\Omega$. The symbol $\mathcal{T}$ will denote regular partitioning of $\Omega$ into elements $K \in \mathcal{T}$, $\bigcup_{K\in\mathcal{T}} \overline{K} = \overline{\Omega}$ and $K \cap K' = \varnothing$ for all $K \neq K' \in \mathcal{T}$. Here, all elements, $K$ will be assumed to have sufficient smooth boundaries $\partial K$, themselves.

For a matrix $\mathsf{N} \in \mathbb{F}^{M\times N}$ of rank $R < \min(M, N)$, the condition number of $\mathsf{N}$ is defined as

$$\mathrm{cond}(\mathsf{N}) = \frac{\sigma_1}{\sigma_R}\,,$$

where $\sigma_1 \geq \sigma_2 \geq \cdots \geq \sigma_R > 0$ are the nonzero singular values of $\mathsf{N}$. As it will be useful later, note that if $\mathsf{N} = \mathsf{M}^*\mathsf{M}$, then $\mathrm{cond}(\mathsf{N}) = \mathrm{cond}(\mathsf{M})^2$ [44].

**Outline.** We attempt to keep the theory of discrete least-squares conspicuously general in Sections 2 and 3. This allows us to draw connections with many other methods in the vast finite element literature as well as highlight the advantages of the DPG setting. In Section 2, we define the Riesz map and introduce notation which allows us to rigorously analyze residual minimization over Hilbert spaces. Here, we give several examples of methods which fall into the DLS framework. In Section 3, we identify DLS with a generalized linear least-squares problem and reduce it to an overdetermined system of equations and associated normal equation. Here, we also demonstrate how to augment these different DLS systems to incorporate boundary conditions as well as compare different direct techniques to solve them. Lastly, we generalize static condensation for the overdetermined system. In Sections 4 and 5, we narrow our focus to DPG methods. In Section 4, we present assembly algorithms in the DPG context. Here, we distinguish between two separate general algorithms; assembly of the normal equation and assembly of the overdetermined system. In Section 5, we illustrate the advantage of the discrete least-squares perspective in several numerical examples including the Poisson equation and the linear acoustics problem. In Section 6 we discuss future possibilities for DLS methods and in Section 7 we summarize our findings.

Appendix A briefly discusses extensions of the DLS approach to constrained minimization problems in quadratic programming and Appendix B provides a second perspective on the DLS static condensation procedure of Section 3.4.

## 2. Minimum residual principles

Due to some technicalities which we wish to avoid for clarity, in this section we consider $\mathbb{F} = \mathbb{R}$. We leave the generalization to $\mathbb{F} = \mathbb{C}$ for the reader. We now begin by defining the vital operators.

### 2.1. **The Riesz map and variational boundary value problems.**

In the abstract setting, a variational boundary value problem is defined with a bilinear form, $b : \mathcal{U} \times \mathcal{V} \to \mathbb{R}$, where $\mathcal{U}$ is the *trial space*, and $\mathcal{V}$ is the *test space*, Assuming that $\mathcal{U}$ and $\mathcal{V}$ are both Hilbert, they are complete in their assigned norms which are induced by inner products and denoted $\|\cdot\|_{\mathcal{U}} = (\cdot,\cdot)_{\mathcal{U}}^{1/2}$ and $\|\cdot\|_{\mathcal{V}} = (\cdot,\cdot)_{\mathcal{V}}^{1/2}$, respectively.

By the Riesz representation theorem, an inner product induces a unique isometric isomorphism between a Hilbert space and its dual [23, 57]. Specifically, letting $\mathcal{W}$ be either $\mathcal{U}$ or $\mathcal{V}$, there exists a unique bijective linear map $\mathcal{R}_{\mathcal{W}} : \mathcal{W} \to \mathcal{W}'$ such that $\| \mathcal{R}_{\mathcal{W}}\, \mathfrak{w}\|_{\mathcal{W}'} = \|\mathfrak{w}\|_{\mathcal{W}}$ for all $\mathfrak{w} \in \mathcal{W}$. This operator is called the Riesz map and is defined by the variational equation

$$\langle \mathcal{R}_{\mathcal{W}}\, \mathfrak{w}, \delta\mathfrak{w}\rangle = (\delta\mathfrak{w}, \mathfrak{w})_{\mathcal{W}} \quad \text{for all } \mathfrak{w}, \delta\mathfrak{w} \in \mathcal{W}\,. \tag{2.1}$$

Following from the theorem, $\|\mathfrak{w}'\|_{\mathcal{W}'} = \| \mathcal{R}_{\mathcal{W}}^{-1}\, \mathfrak{w}'\|_{\mathcal{W}}$ for all $\mathfrak{w}' \in \mathcal{W}'$, $\mathcal{R}_{\mathcal{W}} = \mathcal{R}_{\mathcal{W}}'$, and, moreover, $\mathcal{R}_{\mathcal{W}'} = \mathcal{R}_{\mathcal{W}}^{-1}$ (under the identification $\mathcal{W} \sim \mathcal{W}''$).

Consider, for instance, the circumstance that $\mathcal{W} = L^2(\Omega)$ for some domain $\Omega \subset \mathbb{R}^N$. Here, $\mathcal{R}_{L^2(\Omega)}$ defines the natural identification between $L^2(\Omega)$ and its dual:

$$\langle \mathcal{R}_{L^2(\Omega)}\, v, w\rangle = (w, v)_{L^2(\Omega)} = \int_\Omega wv \quad \text{for all } w, v \in L^2(\Omega)\,. \tag{2.2}$$

Since this is a point of confusion at times, notice that $\mathcal{R}_{L^2(\Omega)}\, v \notin L^2(\Omega)$. That is, it is technically *not* an identity operator, but is instead an element of $(L^2(\Omega))'$.

An abstract variational boundary value problem is a problem of finding a unique trial space element $\mathfrak{u} \in \mathcal{U}$ complementing a given load $\ell \in \mathcal{V}'$ such that

$$b(\mathfrak{u}, \mathfrak{v}) = \ell(\mathfrak{v}) \quad \text{for all } \mathfrak{v} \in \mathcal{V}\,.$$

Similar to the inner product in (2.1), the bilinear form $b$ induces an important linear operator from the trial space to the dual of the test space, $\mathcal{B} : \mathcal{U} \to \mathcal{V}'$:

$$\langle \mathcal{B}\mathfrak{u}, \mathfrak{v}\rangle = b(\mathfrak{u}, \mathfrak{v}) \quad \text{for all } \mathfrak{u} \in \mathcal{U}\,,\ \mathfrak{v} \in \mathcal{V}\,. \tag{2.3}$$

In Section 3, we will see that the Gram matrix and enriched stiffness matrix formed in DLS methods are discrete analogues of Hilbert space operators $\mathcal{R}_{\mathcal{V}}$ and $\mathcal{B}$. First, however, we must explore the idealized minimum residual principles which DLS methods derive from. These principles lead to problems which must be discretized and it is the choice made in their discretization which distinguishes DLS as a unique finite element methodology.

### 2.2. **Minimum residual principles on linear subspaces.**

Here, we derive the idealized Euler-Lagrange equation which DLS emanates from. As we will remark again, later, this framework is presented in a general setting which can be applied to many different problems including that of least-squares finite element methods.

Let $\mathcal{U}_h \subset \mathcal{U}$, be the discrete trial space, where the computed solution will reside. A minimum residual method delivers the optimal discrete solution, $\mathfrak{u}_h^{\mathrm{opt}}$, of the residual minimization problem

$$\mathfrak{u}_h^{\mathrm{opt}} = \underset{\mathfrak{u}_h \in \mathcal{U}_h}{\arg\min}\, \|\mathcal{B}\mathfrak{u}_h - \ell\|_{\mathcal{V}'}^2\,. \tag{2.4}$$

This is the principle supporting Ritz methods, least-squares methods, DPG methods, and general DLS methods, as we will remark later. Using the Riesz map for the test space, (2.4) can be recast in the following form:

$$\mathfrak{u}_h^{\mathrm{opt}} = \operatorname*{arg\,min}_{\mathfrak{u}_h \in \mathcal{U}_h} \left\langle \mathcal{R}_{\mathcal{V}}^{-1}(\mathcal{B}\mathfrak{u}_h - \ell), \mathcal{B}\mathfrak{u}_h - \ell \right\rangle_{\mathcal{V}\times\mathcal{V}'} . \tag{2.5}$$

For the rest of the document we omit the subscripts on the duality pairing $\langle\cdot,\cdot\rangle$, and the spaces are understood based upon the context.

From here, the vanishing Gâteaux derivative at the global minimum implies a variational equation for $\mathfrak{u}_h^{\mathrm{opt}}$, namely,

$$\left\langle \mathcal{R}_{\mathcal{V}}^{-1}\,\mathcal{B}\mathfrak{u}_h^{\mathrm{opt}}, \mathcal{B}\delta\mathfrak{u}_h \right\rangle = \left\langle \mathcal{R}_{\mathcal{V}}^{-1}\,\ell, \mathcal{B}\delta\mathfrak{u}_h \right\rangle \quad \text{for all } \delta\mathfrak{u}_h \in \mathcal{U}_h . \tag{2.6}$$

By introducing the operator $\mathcal{A} = \mathcal{B}'\,\mathcal{R}_{\mathcal{V}}^{-1}\,\mathcal{B}$ and the modified load $f = \mathcal{B}'\,\mathcal{R}_{\mathcal{V}}^{-1}\,\ell \in \mathcal{U}'$, the variational equation (2.6) can be expressed as

$$\left\langle \mathcal{A}\mathfrak{u}_h^{\mathrm{opt}}, \delta\mathfrak{u}_h \right\rangle = f(\delta\mathfrak{u}_h) \quad \text{for all } \delta\mathfrak{u}_h \in \mathcal{U}_h .$$

It is helpful to notice that this new operator, $\mathcal{A} = (\mathcal{R}_{\mathcal{V}}^{-1/2}\,\mathcal{B})'(\mathcal{R}_{\mathcal{V}}^{-1/2}\,\mathcal{B})$, is positive definite and can be identified with an inner product on $\mathcal{U}$,

$$a(\delta\mathfrak{u}, \mathfrak{u}) = \left\langle \mathcal{A}\mathfrak{u}, \delta\mathfrak{u} \right\rangle \quad \text{for all } \mathfrak{u}, \delta\mathfrak{u} \in \mathcal{U} .$$

Therefore, it has the interpretation of a Riesz map on the trial space.

2.3. **Examples of minimum residual principles.** Here, we represent some recognizable methods with the minimum residual principle as well as illustrate some convenient properties of the DPG setting.

*Least-squares.* In a least-squares finite element method, a symmetric bilinear form is usually derived in the following way: Beginning with a linear operator $\mathcal{L} : \mathcal{U} \to L^2(\Omega)$ and a load function $f \in L^2(\Omega)$, seek the solution of the minimization problem

$$\mathfrak{u}_h^{\mathrm{LS}} = \operatorname*{arg\,min}_{\mathfrak{u}_h \in \mathcal{U}_h} \|\mathcal{L}\mathfrak{u}_h - f\|_{L^2(\Omega)}^2 . \tag{2.7}$$

By identifying $\mathcal{V} = L^2(\Omega)$, $\mathcal{B} = \mathcal{R}_{L^2(\Omega)}\,\mathcal{L}$, and $\ell = \mathcal{R}_{L^2(\Omega)}\,f$, we see that (2.4) and (2.7) are equivalent. Indeed, observe that the bilinear form $b$ can be expressed as

$$b(\mathfrak{u}, \mathfrak{v}) = \langle \mathcal{R}_{L^2(\Omega)}\,\mathcal{L}\mathfrak{u}, \mathfrak{v}\rangle = (\mathcal{L}\mathfrak{u}, \mathfrak{v})_{L^2(\Omega)} , \quad \text{for all } \mathfrak{u} \in \mathcal{U} , \ \mathfrak{v} \in \mathcal{V} ,$$

and likewise, $\ell(\mathfrak{v}) = (f, \mathfrak{v})_{L^2(\Omega)}$. Moreover, the residual in (2.4) can be expressed

$$\begin{aligned}
\|\mathcal{B}\mathfrak{u}_h - \ell\|_{L^2(\Omega)'}^2 &= \left\langle \mathcal{R}_{L^2(\Omega)}^{-1}\,\mathcal{R}_{L^2(\Omega)}(\mathcal{L}\mathfrak{u}_h - f), \mathcal{R}_{L^2(\Omega)}(\mathcal{L}\mathfrak{u}_h - f) \right\rangle \\
&= \left\langle \mathcal{R}_{L^2(\Omega)}(\mathcal{L}\mathfrak{u}_h - f), \mathcal{L}\mathfrak{u}_h - f \right\rangle \\
&= (\mathcal{L}\mathfrak{u}_h - f, \mathcal{L}\mathfrak{u}_h - f)_{L^2(\Omega)} \\
&= \|\mathcal{L}\mathfrak{u}_h - f\|_{L^2(\Omega)}^2 .
\end{aligned}$$

Therefore, $\mathfrak{u}_h^{\mathrm{LS}} = \mathfrak{u}_h^{\mathrm{opt}}$. Similarly, variational problem (2.6) is equivalent to

$$(\mathcal{L}\mathfrak{u}_h, \mathcal{L}\delta\mathfrak{u}_h)_{L^2(\Omega)} = (f, \mathcal{L}\delta\mathfrak{u}_h)_{L^2(\Omega)} , \quad \text{for all } \mathfrak{u}_h \in \mathcal{U}_h , \tag{2.8}$$

and the operator $\mathcal{A} = \mathcal{R}_{\mathcal{U}}\,\mathcal{L}^*\mathcal{L}$, where $\mathcal{L}^* = \mathcal{R}_{\mathcal{U}}^{-1}\,\mathcal{L}'\,\mathcal{R}_{L^2(\Omega)} = \mathcal{R}_{\mathcal{U}}^{-1}\,\mathcal{B}'$ is the Hilbert-adjoint operator of $\mathcal{L}$.

*Galerkin methods for symmetric, coercive boundary value problems.* An important variety of linear variational problems, also known as Ritz methods, can be derived from an energy minimization principle

$$\mathfrak{u}_h^{\mathrm{en}} = \operatorname*{arg\,min}_{\mathfrak{u}_h \in \mathcal{U}_h} \mathcal{J}(\mathfrak{u}_h)\,, \tag{2.9}$$

where $\mathcal{J}(\mathfrak{u}) = \frac{1}{2}a(\mathfrak{u},\mathfrak{u}) - \ell(\mathfrak{u})$. Here, $a : \mathcal{U}\times\mathcal{U}\to\mathbb{R}$ is a symmetric, coercive bilinear form, and $\ell \in \mathcal{U}'$ is a linear form. The corresponding variational equation is

$$a(\mathfrak{u}_h^{\mathrm{en}}, \delta\mathfrak{u}_h) = \ell(\delta\mathfrak{u}_h)\,, \quad \text{for all } \delta\mathfrak{u}_h \in \mathcal{U}_h\,. \tag{2.10}$$

In many physical problems, duality methods can be used to equate the minimization principle of (2.9) with a saddle point problem

$$\begin{pmatrix} \mathscr{C}^{-1} & \mathcal{D} \\ \mathcal{D}^* & 0 \end{pmatrix} \begin{pmatrix} \mathfrak{p} \\ \mathfrak{u} \end{pmatrix} = \begin{pmatrix} g \\ f \end{pmatrix}\,,$$

where $\mathcal{D} : \mathcal{U} \to L^2(\Omega)$ is differential operator, $\mathfrak{p} \in L^2(\Omega)$, and $\mathscr{C} : L^2(\Omega) \to L^2(\Omega)$ is symmetric and positive definite. (See [12, 68] for more generality.) In this case,

$$a(\mathfrak{u},\mathfrak{u}) = (\mathscr{C}\mathcal{D}\mathfrak{u}, \mathcal{D}\mathfrak{u})_{L^2(\Omega)} \quad \text{and} \quad \ell(\mathfrak{u}) = (\mathscr{C}\mathcal{D}\mathfrak{u}, g)_{L^2(\Omega)} - (f,\mathfrak{u})_{L^2(\Omega)}\,.$$

If $f = 0$, (2.9) can be identified with the least-squares problem

$$\mathfrak{u}_h^{\mathrm{en}} = \operatorname*{arg\,min}_{\mathfrak{u}_h \in \mathcal{U}_h} \|\mathscr{C}^{1/2}\mathcal{D}\mathfrak{u} - \mathscr{C}^{1/2}g\|^2_{L^2(\Omega)}\,, \tag{2.11}$$

which can, likewise, be equated with a minimum residual principle. Although this formulation has been exploited in some finite element contexts [12, 73], the applicability of an energy principle of the form (2.11) is severely limited.

We now outline a more fruitful procedure. Analogous to (2.3), we may define the operator

$$\langle \mathcal{B}\mathfrak{u}, \mathfrak{v}\rangle = a(\mathfrak{u},\mathfrak{v})\,, \quad \text{for all } \mathfrak{u},\mathfrak{v} \in \mathcal{U}\,.$$

We may also identify $a$ with a norm on the "test space", $\mathcal{V} = \mathcal{U}$, *viz.*,

$$\|\mathfrak{v}\|^2_{\mathcal{V}} = a(\mathfrak{v},\mathfrak{v})\,, \quad \text{for all } \mathfrak{v} \in \mathcal{U}\,.$$

Thus, $\mathcal{R}_{\mathcal{V}} = \mathcal{B}$. With these definitions we see that

$$\|\mathcal{B}\mathfrak{u}_h - \ell\|^2_{\mathcal{V}'} = \langle \mathcal{B}^{-1}(\mathcal{B}\mathfrak{u}_h - \ell), \mathcal{B}\mathfrak{u}_h - \ell\rangle = \langle \mathcal{B}\mathfrak{u}_h - \ell, \mathfrak{u}_h\rangle - \langle \mathcal{B}\mathfrak{u}_h - \ell, \mathcal{B}^{-1}\ell\rangle = 2\mathcal{J}(\mathfrak{u}_h) + \langle \ell, \mathcal{B}^{-1}\ell\rangle.$$

Likewise,

$$\mathfrak{u}_h^{\mathrm{en}} = \operatorname*{arg\,min}_{\mathfrak{u}_h \in \mathcal{U}_h} \mathcal{J}(\mathfrak{u}_h) = \operatorname*{arg\,min}_{\mathfrak{u}_h \in \mathcal{U}_h} \|\mathcal{B}\mathfrak{u}_h - \ell\|^2_{\mathcal{V}'} = \mathfrak{u}_h^{\mathrm{opt}}\,.$$

We have just shown that (2.9) can be expressed as a minimum residual problem and a similar construction is possible with more general (e.g. non-symmetric) variational equations. Unfortunately, the expense of forming a discrete analogue of $\mathcal{R}_{\mathcal{V}}^{-1} = \mathcal{B}^{-1}$ usually surpasses the expense of discretizing (2.10). For this reason, analysis of such a problem in a minimum residual setting is useless unless we move to a broken variational formulation and a DPG method, which we now present the main highlights of.

*DPG methods.* Consider a mesh, $\mathcal{T}$, where $\bigcup_{K\in\mathcal{T}} \overline{K} = \overline{\Omega}$ and $K \cap K' = \emptyset$ for all $K \neq K' \in \mathcal{T}$. We require the following definitions:

- A test space, $\mathcal{V}^{\text{br.}}$, is *broken* if it does not require in its members any form of continuity (conformity) across element boundaries, $\partial K$.
- A test space norm, $\|\cdot\|_{\mathcal{V}^{\text{loc.}}}$, is *localizable* if it induces a finite orthogonal direct sum, $\mathcal{V}^{\text{loc.}} = \bigoplus_{K\in\mathcal{T}} \mathcal{V}_K$, where $\mathfrak{v}_K|_{K'} = 0$ for all $\mathfrak{v}_K \in \mathcal{V}_K$, $K \neq K'$.

In a DPG method, the test space, $\mathcal{V}^{\text{DPG}}$, must be broken and the corresponding test space norm, $\|\cdot\|_{\mathcal{V}^{\text{DPG}}}$, must be localizable. Here, the methodology is often derived with precisely the minimum residual principle (2.4) as well as these two assumptions, although other perspectives are also possible [29].

The broken test space assumption allows the use of discontinuous basis functions in DPG methods, but it is the localizable norm assumption which ensures that the resulting Gram matrix will be block-diagonal (see Section 3.1).

In general, with a slight abuse of notation, the trial space in a DPG method involves two components: $\mathcal{U}^{\text{DPG}} = \mathcal{U}^{\text{fld.}} \times \hat{\mathcal{U}}$, where $\hat{\mathcal{U}}$ can be interpreted as a space of Lagrange multipliers due to the broken nature of the test space. Similarly, the operator $\mathcal{B}^{\text{DPG}}$ in a DPG method can be decomposed into two terms [19]:

$$\mathcal{B}^{\text{DPG}}(\mathfrak{u}^{\text{fld.}}, \hat{\mathfrak{u}}) = \mathcal{B}^{\text{fld.}}\mathfrak{u}^{\text{fld.}} + \hat{\mathcal{B}}\hat{\mathfrak{u}}\,. \tag{2.12}$$

At the infinite-dimensional level, the direct sum decomposition from a localizable norm allows the action of the Riesz map to decouple into its independent element contributions. For instance, given a fixed element $K \in \mathcal{T}$, an element-local test function $\mathfrak{v}_K \in \mathcal{V}_K$, and an arbitrary test function decomposed into element-local contributions, $\delta\mathfrak{v}^{\text{DPG}} = \sum_{K'\in\mathcal{T}} \delta\mathfrak{v}_{K'}$, then

$$\langle \mathcal{R}_{\mathcal{V}^{\text{DPG}}}\mathfrak{v}_K, \delta\mathfrak{v}^{\text{DPG}}\rangle = \sum_{K'\in\mathcal{T}} (\mathfrak{v}_K, \delta\mathfrak{v}_{K'})_{\mathcal{V}^{\text{DPG}}} = (\mathfrak{v}_K, \delta\mathfrak{v}_K)_{\mathcal{V}^{\text{DPG}}} = \langle \mathcal{R}_{\mathcal{V}_K}\mathfrak{v}_K, \delta\mathfrak{v}_K\rangle\,. \tag{2.13}$$

That is, $\mathcal{R}_{\mathcal{V}^{\text{DPG}}}\mathfrak{v}_K : \mathcal{V}^{\text{DPG}} \to \mathbb{F}$ only responds to the nonzero contribution of $\delta\mathfrak{v}^{\text{DPG}}$ within $K$.

Notably, the locality of the Riesz map on a broken test space also exists in its inverse; that is, $\mathcal{R}^{-1}_{\mathcal{V}^{\text{DPG}}} = \bigoplus_{K\in\mathcal{T}} \mathcal{R}^{-1}_{\mathcal{V}_K}$. Therefore, the residual can be decomposed into a *single* sum over the elements of $\mathcal{T}$:

$$\|\mathcal{B}\mathfrak{u}_h - \ell\|^2_{\mathcal{V}^{\text{DPG}\prime}} = \sum_{K\in\mathcal{T}} \langle \mathcal{B}\mathfrak{u}_h - \ell, \mathcal{R}^{-1}_{\mathcal{V}_K}(\mathcal{B}\mathfrak{u}_h - \ell)\rangle = \sum_{K\in\mathcal{T}} \|\mathcal{B}\mathfrak{u}_h - \ell\|^2_{\mathcal{V}'_K}\,. \tag{2.14}$$

Here, each individual term making up the right-hand sum in (2.14) induces an *a posteriori* error estimate which, for each element $K \in \mathcal{T}$, is denoted $\eta_K^2 = \|\mathcal{B}\mathfrak{u}_h - \ell\|^2_{\mathcal{V}'_K}$. These error estimators are often incorporated into adaptive mesh refinement strategies [18].

To identify a connection with least-squares methods, observe that $L^2$ is a broken test space and $\mathcal{R}_{L^2(\Omega)} = \bigoplus_{K\in\mathcal{T}} \mathcal{R}_{L^2(K)}$ and $\hat{\mathcal{B}} = 0$. Likewise,

$$\|\mathcal{L}\mathfrak{u}_h - f\|^2_{L^2(\Omega)} = \sum_{K\in\mathcal{T}} \|\mathcal{L}\mathfrak{u}_h - f\|^2_{L^2(K)}\,,$$

and each localized residual above, $\eta_K = \|\mathcal{L}\mathfrak{u}_h - f\|_{L^2(K)}$, can be used as an *a posteriori* error estimator [6].

2.4. **Non-homogeneous boundary conditions and linear equality constraints.** The minimum residual principle for a variational boundary value problem with non-homogeneous essential boundary conditions generalizes naturally from (2.4). In many physical problems, a selection of components of the trace of the solution are explicitly declared on the boundary of $\Omega$. That is, $\mathfrak{u}_i = \hat{\mathfrak{u}}_i$ on $\partial\Omega$, for some number of prescribed functions $\hat{\mathfrak{u}}_i$, $i \in \mathcal{I} \subset \{1, \ldots, n\}$. The set of all possible prescribed data, in all components, is often a Hilbert space itself—which we will denote $\hat{\mathcal{U}} = \hat{\mathcal{U}}_1 \times \cdots \times \hat{\mathcal{U}}_n$—with an associated continuous linear surjection $\text{tr}_{\mathcal{U}} = \text{tr}_{\mathcal{U}_1} \times \cdots \times \text{tr}_{\mathcal{U}_n} : \mathcal{U} \to \hat{\mathcal{U}}$, called a trace map. With the trace map on $\mathcal{U}$, we may take into account the boundary conditions to hone the solution space into only the affine set $\mathcal{K} = \bigcap_{i\in\mathcal{I}} \text{tr}^{-1}_{\mathcal{U}_i}\{\hat{\mathfrak{u}}_i\} \subset \mathcal{U}$.

Observe that if any single element $\mathfrak{u}^{\text{lift}} \in \mathcal{K}$ is isolated, this so-called "lift" can be used to recharacterize the affine set: $\mathcal{K} = \mathcal{U}^{\text{hom.}} + \mathfrak{u}^{\text{lift}} \subset \mathcal{U}$, where $\mathcal{U}^{\text{hom.}} = \bigcap_{i\in\mathcal{I}} \text{tr}^{-1}_{\mathcal{U}_i}\{0\} = \bigcap_{i\in\mathcal{I}} \text{Null}(\text{tr}_{\mathcal{U}_i})$ is a linear subspace of $\mathcal{U}$. This suggests a procedure which is often performed in precomputation for standard finite element methods. Indeed, it is common practice to isolate a discrete lift, $\mathfrak{u}_h^{\text{lift}}$, which optimally interpolates each component of the solution onto its associated boundary condition in a discrete solution space $\mathcal{U}_h$. (See [27] for an abstract description of traditional finite element methods.) With the lift computed, the problem can then be posed over $\mathcal{K}_h = \mathcal{U}_h^{\text{hom.}} + \mathfrak{u}_h^{\text{lift}} \subset \mathcal{U}_h$, where $\mathcal{U}_h^{\text{hom.}} = \mathcal{U}_h \cap \mathcal{U}^{\text{hom.}}$.

Recall that $b : \mathcal{U} \times \mathcal{V} \to \mathbb{R}$. Abstractly, for such explicit essential boundary conditions, we may therefore begin with an affine subspace $\mathcal{U}_h^{\text{hom.}} + \mathfrak{u}_h^{\text{lift}} \subset \mathcal{U}$ and consider the variational problem

$$\begin{cases} \text{Find } \mathfrak{u}_h \in \mathcal{U}_h^{\text{hom.}} + \mathfrak{u}_h^{\text{lift}} : \\ b(\mathfrak{u}_h, \mathfrak{v}) = \ell(\mathfrak{v}) \quad \text{for all } \mathfrak{v} \in \mathcal{V} . \end{cases} \tag{2.15}$$

The corresponding minimum residual principle is

$$\mathfrak{u}_h^{\text{opt}} = \underset{\mathfrak{u}_h \in \mathcal{U}_h^{\text{hom.}} + \mathfrak{u}_h^{\text{lift}}}{\arg\min} \|\mathcal{B}\mathfrak{u}_h - \ell\|_{\mathcal{V}'}^2 = \underset{\mathfrak{u}_h^{\text{hom.}} \in \mathcal{U}_h^{\text{hom.}}}{\arg\min} \|\mathcal{B}\mathfrak{u}_h^{\text{hom.}} - (\ell - \mathcal{B}\mathfrak{u}_h^{\text{lift}})\|_{\mathcal{V}'}^2 + \mathfrak{u}_h^{\text{lift}} . \tag{2.16}$$

As is plain to see, the minimum residual principle of the explicit boundary condition with an associated lift falls under the subspace minimization theory of Section 2.2.

Unfortunately, not all important boundary conditions can be represented by an explicit prescription of the values of the solution components on $\partial\Omega$. A more general scenario requires the existence of a surjective linear operator $\hat{\mathcal{C}} : \hat{\mathcal{U}} \to \hat{\mathcal{W}}$, where $\hat{\mathcal{W}}$ is an unspecified Hilbert space, for the sake of abstraction. With this operator, we may define a broad range of possible boundary conditions, $\hat{\mathcal{C}}(\mathrm{tr}_{\mathcal{U}}\mathfrak{u}) = \hat{\mathfrak{g}} \in \hat{\mathcal{W}}$. Moreover, we can restrict the possible solutions to the affine subspace $\mathcal{K} = \mathrm{tr}_{\mathcal{U}}^{-1}\{\hat{\mathcal{C}}^{-1}\{\hat{\mathfrak{g}}\}\} \subset \mathcal{U}$.

Now, unless $\hat{\mathcal{C}}$ is the projection map, $\hat{\pi}_{\mathcal{I}} : \hat{\mathcal{U}} \to \prod_{i\in\mathcal{I}} \hat{\mathcal{U}}_i$—in which case we arrive with the explicit boundary conditions above—it is generally difficult to estimate the corresponding lift and homogenous subspace decomposition, $\mathcal{K}_h = \mathcal{U}_h^{\text{hom.}} + \mathfrak{u}_h^{\text{lift}}$, as before. Because of this issue, we require a generalization of the minimum residual principle (2.4). As it extends naturally to the case of more general linear equality constraints that one may wish to satisfy, we deal with the circumstance in full generality.

Consider the linear operator $\mathcal{C} : \mathcal{U} \to \mathcal{W}'$ and the constraint $\mathcal{C}\mathfrak{u} = \mathscr{d}$. A convenient perspective to characterize the associated discrete solution is to search for a unique solution in the implicitly defined, weakly consistent affine subspace $\mathcal{K}_h = \mathcal{K}(\mathcal{U}_h, \mathcal{W}_h) = \{\mathfrak{u}_h \in \mathcal{U}_h \mid \langle \mathcal{C}\mathfrak{u}_h, \mathfrak{w}_h\rangle = \mathscr{d}(\mathfrak{w}_h) \; \forall\, \mathfrak{w}_h \in \mathcal{W}_h\}$. The corresponding discrete variational problem is

$$\begin{cases} \text{Find } \mathfrak{u}_h \in \mathcal{K}_h : \\ b(\mathfrak{u}_h, \mathfrak{v}) = \ell(\mathfrak{v}) \quad \text{for all } \mathfrak{v} \in \mathcal{V} , \end{cases} \tag{2.17}$$

and its complementary constrained minimum residual principle is

$$\mathfrak{u}_h^{\text{opt}} = \underset{\mathfrak{u}_h \in \mathcal{K}_h}{\arg\min} \|\mathcal{B}\mathfrak{u}_h - \ell\|_{\mathcal{V}'}^2 , \tag{2.18}$$

which can equivalently be expressed as

$$\mathfrak{u}_h^{\text{opt}} = \underset{\mathfrak{u}_h \in \mathcal{U}_h}{\arg\min} \sup_{\mathfrak{w}_h \in \mathcal{W}_h} \|\mathcal{B}\mathfrak{u}_h - \ell\|_{\mathcal{V}'}^2 + \langle \mathcal{C}\mathfrak{u}_h - \mathscr{d}, \mathfrak{w}_h\rangle . \tag{2.19}$$

Assuming $\mathcal{K}_h \neq \emptyset$, since $\mathcal{K}_h$ is convex, (2.19) will always have at least one solution [32]. Indeed, existence of a solution to the minimization problem (2.18) can be demonstrated for any convex set $\mathcal{K}_h$ (see Appendix A).

Choosing an inappropriate balance of discrete function spaces $\mathcal{U}_h$ and $\mathcal{W}_h$ can affect the rate of convergence of the optimal solution in (2.19) [12]. In some scenarios, however, $\mathcal{W}$ is finite dimensional and the constraint space $\mathcal{W}_h = \mathcal{W}$ is obvious and will not affect convergence rates when paired with most traditional choices for $\mathcal{U}_h$ (for an example, see [34]).

To avoid potential issues with balancing $\mathcal{U}_h$ and $\mathcal{W}_h$, it is sometimes convenient to consider the following secondary interpretation of optimal solutions to constrained minimization problems. Begin by considering the constraint equation as a second variational problem that we are equally inclined to solve. That is,

$$\begin{cases} \text{Find } \mathfrak{u}_h \in \mathcal{U}_h : \\ b(\mathfrak{u}_h, \mathfrak{v}) = \ell(\mathfrak{v}) \quad \text{for all } \mathfrak{v} \in \mathcal{V} \end{cases} \quad \text{and} \quad \begin{cases} \text{Find } \mathfrak{u}_h \in \mathcal{U}_h : \\ c(\mathfrak{u}_h, \mathfrak{w}) = \mathscr{d}(\mathfrak{w}) \quad \text{for all } \mathfrak{w} \in \mathcal{W} , \end{cases} \tag{2.20}$$

where $c(\mathfrak{u}, \mathfrak{w}) = \langle \mathcal{C}\mathfrak{u}_h, \mathfrak{w}_h \rangle$ is the bilinear form corresponding to the constraint operator. The optimal solution is then the penalized sum of the residuals from both problems:

$$\mathfrak{u}_h^{\text{opt}} = \underset{\mathfrak{u}_h \in \mathcal{U}_h}{\arg\min} \left( \|\mathcal{B}\mathfrak{u}_h - \ell\|^2_{\mathcal{V}'} + \|\mathcal{C}(\mathfrak{u}_h) - d\|^2_{\mathcal{W}'} \right) . \tag{2.21}$$

This is a penalized minimum residual principle. In many scenarios, neither of the problems in (2.20) will have a unique solution. This does not preclude uniqueness of $\mathfrak{u}_h^{\text{opt}}$ in (2.21).

## 3. Discretization

Unlike in the previous section, from now on we consider both cases, $\mathbb{F} = \mathbb{R}$ or $\mathbb{C}$, jointly.

3.1. **The overdetermined system.** Denote an ordered basis for $\mathcal{U}_h$ as $\mathfrak{U}_h = \{\mathfrak{u}_i\}_{i=1}^N$. Similarly, denote a basis for the test space $\mathcal{V}$ as $\mathfrak{V} = \{\mathfrak{v}_\imath\}_{\imath \in \mathcal{I}}$ and denote the first $M > N$ of a countable subset of these functions as the ordered basis, $\mathfrak{V}_r = \{\mathfrak{v}_i\}_{i=1}^M$, for the discrete test space $\mathcal{V}_r$. Likewise, the discrete Riesz map, $\mathcal{R}_{\mathcal{V}_r} : \mathcal{V}_r \to \mathcal{V}_r'$, is defined

$$\langle \mathcal{R}_{\mathcal{V}_r} \mathfrak{v}, \delta\mathfrak{v} \rangle = (\delta\mathfrak{v}, \mathfrak{v})_{\mathcal{V}} \quad \text{for all } \mathfrak{v}, \delta\mathfrak{v} \in \mathcal{V}_r = \operatorname{span}\left(\{\mathfrak{v}_i\}_{i=1}^M\right) .$$

With this definition, we formulate the practical minimization problem

$$\mathfrak{u}_{h,r}^{\text{opt}} = \underset{\mathfrak{u}_h \in \mathcal{U}_h}{\arg\min} \left\| \mathcal{R}_{\mathcal{V}_r}^{-1}(\mathcal{B}\mathfrak{u}_h - \ell) \right\|^2_{\mathcal{V}} . \tag{3.1}$$

In (3.1) above, $\mathcal{R}_{\mathcal{V}_r}^{-1} : \mathcal{V}' \to \mathcal{V}$ is an abuse of notation naturally representing the extension of $\mathcal{R}_{\mathcal{V}_r}^{-1} : \mathcal{V}_r' \to \mathcal{V}_r$ by zero. Or, in other terms, $\mathcal{R}_{\mathcal{V}_r}^{-1} = \mathcal{R}_{\mathcal{V}}^{-1} \circ \iota_r'$, where $\iota_r'$ is the dual of the canonical embedding, $\iota_r : \mathcal{V}_r \to \mathcal{V}$.

Although (3.1) generally represents the quadratic minimization problem actually solved in a computer implementation, it is perhaps more naturally expressed in terms of the actual finite element matrices. In such a representation, it can be identified with a (discrete) generalized least-squares problem on the coefficients, $\mathsf{u} = [\mathsf{u}_i]_{i=1}^N \in \mathbb{F}^N$, of the solution $\mathfrak{u}_h \in \mathcal{U}_h$ represented in the $\{\mathfrak{u}_i\}_{i=1}^N$ basis,

$$\mathfrak{u}_h = \sum_{i=1}^N \mathsf{u}_i \mathfrak{u}_i \in \mathcal{U}_h .$$

Indeed, define $\mathsf{B}_{ij} = b(\mathfrak{u}_j, \mathfrak{v}_i)$, $\mathsf{G}_{ij} = (\mathfrak{v}_i, \mathfrak{v}_j)_{\mathcal{V}}$, and $\mathsf{l}_i = \ell(\mathfrak{v}_i)$ for each ordered basis element $\mathfrak{u}_j \in \mathfrak{U}$ and $\mathfrak{v}_i, \mathfrak{v}_j \in \mathfrak{V}_r$, and define the optimal coefficients to be

$$\mathsf{u}^{\text{opt}} = \underset{\mathsf{u} \in \mathbb{F}^N}{\arg\min}\, (\mathsf{B}\mathsf{u} - \mathsf{l})^* \mathsf{G}^{-1} (\mathsf{B}\mathsf{u} - \mathsf{l}) . \tag{3.2}$$

Of course, for $\mathsf{u}^{\text{opt}} = [\mathsf{u}_i^{\text{opt}}]_{i=1}^N$,

$$\mathfrak{u}_{h,r}^{\text{opt}} = \sum_{i=1}^N \mathsf{u}_i^{\text{opt}} \mathfrak{u}_i ,$$

as the reader may verify.

Now, reflect upon the matrix $\mathsf{G}$. In place of a general basis $\mathfrak{V}_r$, had we used an orthonormal basis $\widetilde{\mathfrak{V}}_r = \{\widetilde{\mathfrak{v}}_i\}_{i=1}^M$, then minimization problem (3.2) would be equivalent to an ordinary linear least-squares problem. That is,

$$(\widetilde{\mathsf{B}}\mathsf{u} - \widetilde{\mathsf{l}})^* \widetilde{\mathsf{G}}^{-1} (\widetilde{\mathsf{B}}\mathsf{u} - \widetilde{\mathsf{l}}) = \|\widetilde{\mathsf{B}}\mathsf{u} - \widetilde{\mathsf{l}}\|_2^2 ,$$

where $\widetilde{\mathsf{B}}_{ij} = b(\mathfrak{u}_j, \widetilde{\mathfrak{v}}_i)$, $\widetilde{\mathsf{G}}_{ij} = (\widetilde{\mathfrak{v}}_i, \widetilde{\mathfrak{v}}_j)_{\mathcal{V}} = \delta_{ij}$, and $\widetilde{\mathsf{l}}_i = \ell(\widetilde{\mathfrak{v}}_i)$, and therefore

$$\mathsf{u}^{\text{opt}} = \underset{\mathsf{u} \in \mathbb{F}^N}{\arg\min}\, \|\widetilde{\mathsf{B}}\mathsf{u} - \widetilde{\mathsf{l}}\|_2^2 . \tag{3.3}$$

As we will now demonstrate, precomputing an orthonormal basis $\widetilde{\mathfrak{V}}_r$ is not necessary and can be performed implicitly in the process of inverting $\mathsf{G}$.

Indeed, since $\mathsf{G}$ is positive definite, it is congruent to the identity matrix, so one can arrive at (3.3) by defining $\widetilde{\mathsf{B}} = \mathsf{W}^{-1}\mathsf{B}$ and $\widetilde{\mathsf{l}} = \mathsf{W}^{-1}\mathsf{l}$, where $\mathsf{W}$ is any matrix solving the equation $\mathsf{W}\mathsf{W}^* = \mathsf{G}$. Such a factorization is naturally related to finding a discretely $\mathcal{V}$-orthogonal change-of-basis transformation. For instance, consider an arbitrary basis vector $\mathfrak{v}_i \in \mathfrak{V}_r$ represented in an orthonormal basis $\widetilde{\mathfrak{V}}_r$,

$$\mathfrak{v}_i = \sum_{j=1}^{M} \mathsf{W}_{ij} \widetilde{\mathfrak{v}}_j \,,$$

where $(\widetilde{\mathfrak{v}}_i, \widetilde{\mathfrak{v}}_j)_{\mathcal{V}} = \delta_{ij}$ for all $\widetilde{\mathfrak{v}}_i, \widetilde{\mathfrak{v}}_j \in \widetilde{\mathfrak{V}}_r$. With this expression, the Gram matrix can be represented as

$$\mathsf{G}_{ij} = (\mathfrak{v}_i, \mathfrak{v}_j)_{\mathcal{V}} = \sum_{k,\,l=1}^{M} \mathsf{W}_{ik} \overline{\mathsf{W}}_{jl} (\widetilde{\mathfrak{v}}_k, \widetilde{\mathfrak{v}}_l) = \sum_{k=1}^{M} \mathsf{W}_{ik} \overline{\mathsf{W}}_{jk} \,.$$

i.e. $\mathsf{G} = \mathsf{W}\mathsf{W}^*$.

If $M > 1$, the equation $\mathsf{W}\mathsf{W}^* = \mathsf{G}$ has an infinite number of solutions but a convenient choice is the (unique) Cholesky factorization, $\mathsf{G} = \mathsf{L}\mathsf{L}^*$. For lower-triangular matrix $\mathsf{W}_{\text{Chol.}} = \mathsf{L}$, we can rewrite

$$\begin{aligned} \mathbf{u}^{\text{opt}} &= \arg\min_{\mathbf{u}\in\mathbb{F}^N} \big(\mathsf{W}_{\text{Chol.}}^{-1}(\mathsf{B}\mathbf{u} - \mathsf{l})\big)^* \mathsf{W}_{\text{Chol.}}^{-1}(\mathsf{B}\mathbf{u} - \mathsf{l}) \\ &= \arg\min_{\mathbf{u}\in\mathbb{F}^N} \big\|\mathsf{L}^{-1}(\mathsf{B}\mathbf{u} - \mathsf{l})\big\|_2^2 \,, \end{aligned} \tag{3.4}$$

and the matrices $\widetilde{\mathsf{B}} = \mathsf{L}^{-1}\mathsf{B}$ and $\widetilde{\mathsf{l}} = \mathsf{L}^{-1}\mathsf{l}$ can be efficiently computed by back-substitution. This is choice is usually assumed throughout this paper. In the statistics and signal processing communities, the factoring and row-weighting procedure described above is known as “whitening” or “sphering” and we refer the reader to [50] to explore the benefits of the other possibilities for the change-of-basis matrices $\mathsf{W}$ in that context.

*Synopsis.* Ultimately, discretization of the Riesz map in a minimum residual principle induces a generalized least-squares problem

$$\min_{\mathbf{u}\in\mathbb{F}^N} \|\mathsf{B}\mathbf{u} - \mathsf{l}\|_{\mathsf{G}^{-1}} \,,$$

where $\mathsf{G}$, the Gram matrix, is the discretization of the Riesz map. This is equivalent to the special equilibrium system, or saddle-point problem,

$$\begin{bmatrix} \mathsf{G} & \mathsf{B} \\ \mathsf{B}^* & 0 \end{bmatrix} \begin{bmatrix} \mathbf{p} \\ \mathbf{u} \end{bmatrix} = \begin{bmatrix} \mathsf{l} \\ 0 \end{bmatrix} . \tag{3.5}$$

Since $\mathbf{u}$ is the subvector of primary importance, the statically-condensed equations

$$\mathsf{A}\mathbf{u} = \mathsf{f} \,, \tag{3.6}$$

where $\mathsf{A} = \mathsf{B}^*\mathsf{G}^{-1}\mathsf{B}$ and $\mathsf{f} = \mathsf{B}^*\mathsf{G}^{-1}\mathsf{l}$, can be used to characterize the solution. If $\mathsf{G}$ can be efficiently factored $\mathsf{G} = \mathsf{W}\mathsf{W}^*$—this is the case which we are predominantly interested in—it is also useful to rewrite this problem as

$$\min_{\mathbf{u}\in\mathbb{F}^N} \big\|\mathsf{W}^{-1}(\mathsf{B}\mathbf{u} - \mathsf{l})\big\|_2^2 \,.$$

The normal equation of are then (3.6).

We will discuss efficient solution algorithms for these equations in Section 3.3, but first we mention some special properties of the DPG setting for this class of problems and also discuss the imposition of boundary conditions and linear constraints in general terms.

*DLS for traditional methods.* Beginning with any standard finite element method, it is tempting to consider the corresponding stiffness matrix $\mathsf{B} \in \mathbb{F}^{M\times N}$ and load vector $\mathsf{l} \in \mathbb{F}^N$, with $M = N$. Defining the Gram matrix to be the identity matrix, $\mathsf{G} = \mathrm{I}$—or any invertible matrix, really—the full rank linear system of equations

$$\mathsf{B}\mathbf{u} = \mathsf{l} \,,$$

is equivalent to each of the formulations above.

In some methods, it is also quite straightforward to extend the associated test space. Therefore, let us assume again that $M > N$. In this situation, the naive choice $\mathsf{G} = \mathrm{I}$ naturally induces the linear system (3.6), but the condition number of $\mathsf{A} = \mathsf{B}^*\mathsf{B}$ would then be the square of the condition number of $\mathsf{B}$. Since the condition number of $\mathsf{B}$ is often $\mathcal{O}(h^{-2})$, this would usually be a gross impediment to practicality. Fortunately, with an appropriate choice of test norm or Gram matrix, the condition number can be far better behaved. Indeed, in Section 5 we demonstrate condition number growth of $\mathsf{A}$ to be $\mathcal{O}(h^{-2}) \neq \mathcal{O}(h^{-4})$—and, thus, $\mathrm{cond}(\widetilde{B}) = \mathcal{O}(h^{-1})$—in several practical settings.

*DPG methods.* As mentioned previously, if the test space is broken and localizable (see Section 2.3), the Gram matrix $\mathsf{G}$ can be constructed to be block diagonal. Indeed, let $\mathcal{V}^{\mathrm{DPG}} = \bigoplus_{K\in\mathcal{T}} \mathcal{V}_K$ be the DPG test space and let $\mathfrak{V}_r^{\mathrm{DPG}} = \bigoplus_{K\in\mathcal{T}} \mathfrak{V}_K$, $\mathfrak{V}_K \subset \mathcal{V}_K$, be a broken and localized basis for the finite-dimensional subspace $\mathcal{V}_r^{\mathrm{DPG}} \subset \mathcal{V}^{\mathrm{DPG}}$. That is, for each $i \in \{1,\cdots,M\}$, $\mathfrak{v}_i \in \mathcal{V}_K$ for some $K \in \mathcal{T}$ where $\mathfrak{v}_i|_{K'} = 0$ for all $K' \neq K$.

Fixing the previous $i$ and $K$, observe that, for each $j \in \{1,\cdots,M\}$,

$$\mathsf{G}_{ij} = (\mathfrak{v}_i, \mathfrak{v}_j)_{\mathcal{V}^{\mathrm{DPG}}} = (\mathfrak{v}_i, \mathfrak{v}_j)_{\mathcal{V}_K}\,,$$

by (2.13). Therefore, $\mathsf{G} = \mathrm{diag}(\mathsf{G}_K)$, where each $\mathsf{G}_K$ is the local Gram matrix,

$$(\mathsf{G}_K)_{ij} = (\mathfrak{v}_i, \mathfrak{v}_j)_{\mathcal{V}_K}\,,$$

for each $\mathfrak{v}_i, \mathfrak{v}_j \in \mathfrak{V}_K$.

As also alluded to previously, in some circumstances, many, or all, $\mathsf{G}_K$ can be made be *diagonal* if the basis is properly compatible with the test norm. For regular enough elements, this is often possible when $\mathcal{V} = L^2$—indeed, this was verified in the FOSLS experiments in Section 5.1. The interested reader may consult [7, 36] and references therein, to explore such orthogonal or sparsity-optimized shape functions.

Recall from (2.12), $\mathcal{B}^{\mathrm{DPG}}(\mathfrak{u}^{\mathrm{fld.}}, \hat{\mathfrak{u}}) = \mathcal{B}^{\mathrm{fld.}}\mathfrak{u}^{\mathrm{fld.}} + \hat{\mathcal{B}}\hat{\mathfrak{u}}$. Through an abuse of notation, denote the coefficients of $\mathfrak{u}^{\mathrm{fld.}}$ as $\mathbf{u}$ and the coefficients of $\hat{\mathfrak{u}}$ as $\hat{\mathbf{u}}$. Similarly define $\mathsf{B}$ and $\hat{\mathsf{B}}$. With these definitions, the DPG analog of the saddle-point problem (3.5) is

$$\begin{bmatrix} \mathsf{G} & \mathsf{B} & \hat{\mathsf{B}} \\ \mathsf{B}^* & 0 & 0 \\ \hat{\mathsf{B}}^* & 0 & 0 \end{bmatrix} \begin{bmatrix} \mathbf{p} \\ \mathbf{u} \\ \hat{\mathbf{u}} \end{bmatrix} = \begin{bmatrix} \mathsf{l} \\ 0 \\ 0 \end{bmatrix}. \tag{3.7}$$

3.2. **Boundary conditions and linear equality constraints.** Recall the matrices $\mathsf{G}$, $\mathsf{B}$, $\widetilde{\mathsf{B}}$, and $\mathsf{A}$, and vectors $\mathsf{l}$, $\widetilde{\mathsf{l}}$, and $\mathsf{f}$ from Section 3.1. As demonstrated in Section 2.4, in principle, a constrained minimum residual problem will involve two test spaces, $\mathcal{V}$ and $\mathcal{W}$. We have already fixed a truncated basis for $\mathcal{V}$ spanning $\mathcal{V}_r \subset \mathcal{V}$. Here, we additionally require a similar basis for $\mathcal{W}_s \subset \mathcal{W}$, denoted $\mathfrak{W}_s = \{\mathfrak{w}_i\}_{i=1}^L$. Drawing from $\mathfrak{W}_s$, we now define $\mathsf{C}_{ij} = c(\mathfrak{w}_j, \mathfrak{w}_i)$, $\mathsf{H}_{ij} = (\mathfrak{w}_i, \mathfrak{w}_j)_{\mathcal{W}}$, and $\mathsf{d}_i = d(\mathfrak{w}_i)$.

Only in discretizing the penalized problem, (2.21), will we require the second Gram matrix, $\mathsf{H}$. In the saddle-point problem, (2.19), however, it will be required to solve for the constraint coefficients, $\mathbf{w} = [\mathsf{w}_i]_{i=1}^L \in \mathbb{F}^L$.

Proceeding as before, the saddle-point problem can be expressed

$$\mathbf{u}^{\mathrm{opt}} = \underset{\mathbf{u}\in\mathbb{F}^N}{\arg\min} \sup_{\mathbf{w}\in\mathbb{F}^L} \left( \|\mathsf{B}\mathbf{u} - \mathsf{l}\|_{\mathsf{G}^{-1}} + \mathbf{w}^*(\mathsf{C}\mathbf{u} - \mathsf{d}) \right). \tag{3.8}$$

Or, in block matrix form, as

$$\begin{bmatrix} \mathsf{G} & \mathsf{B} & 0 \\ \mathsf{B}^* & 0 & \mathsf{C}^* \\ 0 & \mathsf{C} & 0 \end{bmatrix} \begin{bmatrix} \mathbf{p}^{\mathrm{opt}} \\ \mathbf{u}^{\mathrm{opt}} \\ \mathbf{w}^{\mathrm{opt}} \end{bmatrix} = \begin{bmatrix} \mathsf{l} \\ 0 \\ \mathsf{d} \end{bmatrix} \quad \Longleftrightarrow \quad \begin{bmatrix} \mathsf{A} & \mathsf{C}^* \\ \mathsf{C} & 0 \end{bmatrix} \begin{bmatrix} \mathbf{u}^{\mathrm{opt}} \\ \mathbf{w}^{\mathrm{opt}} \end{bmatrix} = \begin{bmatrix} \mathsf{f} \\ \mathsf{d} \end{bmatrix}. \tag{3.9}$$

Notably, this problem is well-posed if and only if the linear system $\mathsf{C}\mathbf{u} = \mathsf{d}$ is consistent, $\mathsf{C}$ has full row rank, and $\mathrm{Null}(\mathsf{A}) \cap \mathrm{Null}(\mathsf{C}) = \{0\}$ [8, Chapter 5]. Moreover, the right-hand equation in (3.9) is fundamentally different from the DLS saddle-point problem (3.5), even when $\mathsf{d} = 0$. This is not only because we cannot expect $\mathsf{A}$ to have any easily exploitable structure, but because we are predominately interested in the first subvector, $\mathbf{u}^{\mathrm{opt}}$, not $\mathbf{w}^{\mathrm{opt}}$.

Alternatively, the penalized problem is

$$\mathbf{u}^{\mathrm{opt}} = \underset{\mathbf{u}\in\mathbb{F}^N}{\arg\min} \left( \|\mathsf{B}\mathbf{u} - \mathsf{l}\|^2_{\mathsf{G}^{-1}} + \|\mathsf{C}\mathbf{u} - \mathsf{d}\|^2_{\mathsf{H}^{-1}} \right), \tag{3.10}$$

resulting in the saddle point problem

$$\begin{bmatrix} \mathsf{G} & \mathsf{B} & 0 \\ \mathsf{B}^* & 0 & \mathsf{C}^* \\ 0 & \mathsf{C} & \mathsf{H} \end{bmatrix} \begin{bmatrix} \mathbf{p}^{\mathrm{opt}} \\ \mathbf{u}^{\mathrm{opt}} \\ \mathbf{w}^{\mathrm{opt}} \end{bmatrix} = \begin{bmatrix} \mathsf{l} \\ 0 \\ \mathsf{d} \end{bmatrix}, \tag{3.11}$$

or, equivalently,

$$(\mathsf{A} + \mathsf{E})\mathbf{u}^{\mathrm{opt}} = \mathsf{f} + \mathsf{e}\,, \tag{3.12}$$

where $\mathsf{E} = \mathsf{C}^*\mathsf{H}^{-1}\mathsf{C}$ and $\mathsf{e} = \mathsf{C}^*\mathsf{H}^{-1}\mathsf{d}$. This is well-posed if and only if $\mathrm{Null}(\mathsf{A}) \cap \mathrm{Null}(\mathsf{C}) = \{0\}$.

Most often, the Gram matrix $\mathsf{H}$ has a tunable penalization parameter: $\mathsf{H} = \alpha^{-2}\mathsf{H}_0$. The magnitude of $\alpha$ determines the bias on the constraint in (3.10). Intuitively, the solution of (3.11) converges to the solution of (3.9) as $\alpha \to \infty$, and indeed, this convergence can be made precise [5]. Unfortunately, as $\alpha$ grows, so will the condition number of $\mathsf{E}$, and therefore, so will the round-off error in (3.12).

Often, $\mathsf{H}_0 = \mathrm{I}$. In this scenario, (3.11) can also be expressed as the least-squares problem for an overdetermined system:

$$\min_{\mathbf{u}\in\mathbb{F}^N} \left\| \begin{bmatrix} \alpha\mathsf{C} \\ \widetilde{\mathsf{B}} \end{bmatrix} \mathbf{u} - \begin{bmatrix} \alpha\mathsf{d} \\ \widetilde{\mathsf{l}} \end{bmatrix} \right\|_2, \tag{3.13}$$

which is less sensitive to the penalization parameter than the normal equation when solved with QR. In the least-squares context, this is referred to as the "method of weighting" [71]. And, the extension to scenarios when the Gram matrix can be factored, $\mathsf{H} = \mathsf{M}\mathsf{M}^*$, should be clear.

In the special, but common, case of equality constraints in the solution coefficients, we can decompose the vector $\mathbf{u} = [\mathbf{u}_{\mathrm{hom.}} \mid \mathbf{u}_{\mathrm{lift}}]^{\mathsf{T}}$ into the *free* homogeneous coefficients, $\mathbf{u}_{\mathrm{hom.}}$, and the *fixed* lifted coefficients, $\mathbf{u}_{\mathrm{lift}}$. We may then decompose the identity matrix $\mathrm{I} = [\mathsf{Q}_{\mathrm{hom.}} \mid \mathsf{Q}_{\mathrm{lift}}]$ into two submatrices with (clearly) orthonormal columns, $\mathsf{Q}^*_{\mathrm{hom.}}\mathbf{u} = \mathbf{u}_{\mathrm{hom.}}$ and $\mathsf{Q}^*_{\mathrm{lift}}\mathbf{u} = \mathbf{u}_{\mathrm{lift}}$.

Now, observe that

$$\mathsf{B}\mathbf{u} = \mathsf{B}\mathsf{Q}_{\mathrm{hom.}}\mathsf{Q}^*_{\mathrm{hom.}}\mathbf{u} + \mathsf{B}\mathsf{Q}_{\mathrm{lift}}\mathsf{Q}^*_{\mathrm{lift}}\mathbf{u} = \mathsf{B}_{\mathrm{hom.}}\mathbf{u}_{\mathrm{hom.}} + \mathsf{B}_{\mathrm{lift}}\mathbf{u}_{\mathrm{lift}}\,,$$

where $\mathsf{B}_{\mathrm{hom.}} = \mathsf{B}\mathsf{Q}_{\mathrm{hom.}}$ and $\mathsf{B}_{\mathrm{lift}} = \mathsf{B}\mathsf{Q}_{\mathrm{lift}}$. Then, the optimal solution coefficients can then be characterized as

$$\mathbf{u}^{\mathrm{opt}} = \mathop{\arg\min}_{\mathbf{u}_{\mathrm{hom.}}\in\mathbb{F}^{N_{\mathrm{hom.}}}} \|\mathsf{B}_{\mathrm{hom.}}\mathbf{u}_{\mathrm{hom.}} - \mathsf{l}_{\mathrm{lift}}\|_{\mathsf{G}^{-1}} + \mathbf{u}_{\mathrm{lift}} = \mathop{\arg\min}_{\mathbf{u}_{\mathrm{hom.}}\in\mathbb{F}^{N_{\mathrm{hom.}}}} \|\widetilde{\mathsf{B}}_{\mathrm{hom.}}\mathbf{u}_{\mathrm{hom.}} - \widetilde{\mathsf{l}}_{\mathrm{lift}}\|_2 + \mathbf{u}_{\mathrm{lift}}\,,$$

where $\mathsf{l}_{\mathrm{lift}} = \mathsf{l} - \mathsf{B}_{\mathrm{lift}}\mathbf{u}_{\mathrm{lift}}$, $\widetilde{\mathsf{B}}_{\mathrm{hom.}} = \widetilde{\mathsf{B}}\mathsf{Q}_{\mathrm{hom.}}$, and $\widetilde{\mathsf{l}}_{\mathrm{lift}} = \widetilde{\mathsf{l}} - \widetilde{\mathsf{B}}\mathsf{Q}_{\mathrm{lift}}\mathbf{u}_{\mathrm{lift}}$.

In the construction above, introducing the matrices $\mathsf{Q}_{\mathrm{hom.}}$ and $\mathsf{Q}_{\mathrm{lift}}$ was seemingly unnecessary because $\mathsf{B} = [\mathsf{B}_{\mathrm{hom.}} \mid \mathsf{B}_{\mathrm{lift}}]$. However, the solution coefficients will generally not be ordered as above and this detail allows us to connect to the more general constraint case, by identifying $\mathsf{C}$ with $\mathsf{Q}^*_{\mathrm{lift}}$:

$$\begin{bmatrix} \mathsf{G} & \mathsf{B} & 0 \\ \mathsf{B}^* & 0 & \mathsf{Q}_{\mathrm{lift}} \\ 0 & \mathsf{Q}^*_{\mathrm{lift}} & 0 \end{bmatrix} \begin{bmatrix} \mathbf{p}^{\mathrm{opt}} \\ \mathbf{u}^{\mathrm{opt}} \\ \mathbf{w}^{\mathrm{opt}} \end{bmatrix} = \begin{bmatrix} \mathsf{l} \\ 0 \\ \mathbf{u}^{\mathrm{lift}} \end{bmatrix} \quad \Longleftrightarrow \quad \begin{bmatrix} \mathsf{A} & \mathsf{Q}_{\mathrm{lift}} \\ \mathsf{Q}^*_{\mathrm{lift}} & 0 \end{bmatrix} \begin{bmatrix} \mathbf{u}^{\mathrm{opt}} \\ \mathbf{w}^{\mathrm{opt}} \end{bmatrix} = \begin{bmatrix} \mathsf{f} \\ \mathbf{u}^{\mathrm{lift}} \end{bmatrix}.$$

After multiplying the first row of the second equation by $\mathsf{Q}^*_{\mathrm{hom.}}$, it is easy to see that these equations imply

$$\mathsf{A}_{\mathrm{hom.}}\mathbf{u}^{\mathrm{opt}}_{\mathrm{hom.}} = \mathsf{B}^*_{\mathrm{hom.}}\mathsf{G}^{-1}\mathsf{l}_{\mathrm{lift}}\,,$$

where $\mathsf{A}_{\mathrm{hom.}} = \mathsf{Q}^*_{\mathrm{hom.}}\mathsf{A}\mathsf{Q}_{\mathrm{hom.}} = \mathsf{B}^*_{\mathrm{hom.}}\mathsf{G}^{-1}\mathsf{B}_{\mathrm{hom.}}$.

*Post-processing precomputed solutions.* In some situations, large subvectors of the coefficients of the solution, $\mathbf{u}_{\mathrm{sol.}}$, corresponding to entire components, $\mathfrak{u}_i$, $i \in \mathcal{I}_{\mathrm{sol.}} \subsetneq \{1, \ldots, n\}$, of the discrete solution may already be available at the outset of computation. For instance, this is possible if an auxiliary finite element method was used to compute $\mathfrak{u}_i$. In this scenario, the remaining solution coefficients can be computed using the minimum residual framework with equality constraints as described above and defining the associated $\mathsf{Q}_{\mathrm{sol.}} = \mathsf{C}^*$, $\mathsf{C}\mathbf{u} = \mathbf{u}_{\mathrm{sol.}}$.

3.3. **Solution algorithms.** There are several prevalent strategies to solve generalized linear least-squares problems

$$\text{(DLS)} \qquad \min_{\mathbf{u}\in\mathbb{F}^N} \|\mathsf{B}\mathbf{u} - \mathsf{l}\|_{\mathsf{G}^{-1}} .$$

Under infinite precision arithmetic, they are each essentially equivalent, however, algorithmically speaking, each approach is significantly different. We briefly survey the most important direct methods in this subsection.

*The normal equation.* In order to solve for $\mathbf{u}^{\mathrm{opt}}$, we could begin by forming the normal equation:

$$\text{(NE)} \qquad \mathsf{A}\mathbf{u}^{\mathrm{opt}} = \mathsf{f}\,, \quad \text{where} \quad \mathsf{A} = \mathsf{B}^*\mathsf{G}^{-1}\mathsf{B} \quad \text{and} \quad \mathsf{f} = \mathsf{B}^*\mathsf{G}^{-1}\mathsf{l}\,.$$

Beginning in this way is standard practice in the DPG community and is advantageous in that the stiffness matrix and load vector, $\mathsf{A}$ and $\mathsf{f}$, can be constructed locally and assembled in sparse format using standard finite element assembly algorithms. Moreover, $\mathsf{A} = \widetilde{\mathsf{B}}^*\widetilde{\mathsf{B}}$ is Hermitian positive-definite and therefore has a structure amenable to efficient linear solvers not usually available for many challenging problems. Unfortunately, the condition number of $\mathsf{A}$ will grow quadratically with the condition number of $\widetilde{\mathsf{B}}$ and, likewise, so will the upper bound on the round-off error of the normal equation solution. However, this is not to say that when the normal equation are formed the growth of the condition number will be expected to be faster than most standard finite element methods. For instance, when using a first-order system (FOS) formulation, some methods using the normal equation can be proven to have condition number growth $\mathcal{O}(h^{-2})$, where $h$ is the element size in a quasi-uniformly refined mesh [45]. Similarly, this has also been shown to be true in a more general context for the linear systems arising from FOS least-squares (FOSLS) methods [10, 11] which the normal equation above are related to (see Section 5.1).

Nevertheless, the scaling constant controlling the condition number of the stiffness matrix is often large in a FOS setting due to the additional equations and unknowns. This is an impediment to producing accurate solutions to very difficult problems and can induce large and sometimes overwhelming numerical errors. This is one reason that it is sometimes convenient to consider other approaches which deal explicitly with the matrices $\mathsf{B}$, $\mathsf{W}$, and $\mathsf{l}$, and avoid the normal equation altogether.

*Orthogonal decomposition methods.* The most practical alternative to the normal equation when solving for $\mathbf{u}^{\mathrm{opt}}$ is to deal directly with the matrices $\widetilde{\mathsf{B}}$ and $\widetilde{\mathsf{l}}$ coming from the (sparsely weighted) linear least-squares problem

$$\text{(LS)} \qquad \min_{\mathbf{u}\in\mathbb{F}^N} \|\widetilde{\mathsf{B}}\mathbf{u} - \widetilde{\mathsf{l}}\|_2\,, \quad \text{where} \quad \widetilde{\mathsf{B}} = \mathsf{W}^{-1}\mathsf{B}\,, \quad \widetilde{\mathsf{l}} = \mathsf{W}^{-1}\mathsf{l}\,, \quad \text{and} \quad \mathsf{W}\mathsf{W}^* = \mathsf{G}\,.$$

The most common of these approaches is the orthogonalization algorithm called QR-factorization (Householder, Givens, MGS) first introduced for least-squares problems in [41]. Other direct approaches are SVD, complete orthogonal decomposition, and Peter-Wilkinson as well as various hybrid methods [8]. Each of these approaches are usually less efficient than solving via the normal equation, but are often preferred because they are more numerically stable.

For the purposes we have in mind, the matrix $\widetilde{\mathsf{B}}$ will be large and sparse and so, because not all of the methods above are well suited for sparse matrices or amenable to parallel computing, we will focus only on the QR approach. As shown in various textbooks [8, 44, 70], the relative error in the solution from a least-squares QR solve is controlled by $\mathrm{cond}(\widetilde{\mathsf{B}}) + \rho(\widetilde{\mathsf{B}},\widetilde{\mathsf{l}}) \cdot \mathrm{cond}(\widetilde{\mathsf{B}})^2$, where

$$\rho(\widetilde{\mathsf{B}},\widetilde{\mathsf{l}}) = \frac{\|\widetilde{\mathsf{B}}\mathbf{u}^{\mathrm{opt}} - \widetilde{\mathsf{l}}\|_2}{\|\widetilde{\mathsf{B}}\|_2\|\mathbf{u}^{\mathrm{opt}}\|_2}\,, \tag{3.14}$$

and therefore depends upon the load vector.

Due to consistency of the problem, $\ell \in \mathrm{Range}(\mathcal{B})$. Therefore, the residual (the numerator in (3.14)) is expected to tend to zero as the mesh is refined. That is, $\rho(\widetilde{\mathsf{B}},\widetilde{\mathsf{l}}) \to 0$ as $h \to 0$. Indeed, the reader may observe that the $\rho(\widetilde{\mathsf{B}},\widetilde{\mathsf{l}})$ even vanishes if $\widetilde{\mathsf{l}} \in \mathrm{Range}(\widetilde{\mathsf{B}})$. Of course, the validity of this convergence as well as its rate will be determined by the interpolation spaces used in the discretization. However, in many common scenarios the *a priori* bound can be proven to decrease at a rate of at least $\mathcal{O}(h)$. Indeed, for many cases we have in mind, this rate is of the form $\mathcal{O}(h^p)$ where $p \geq 1$ is a parameter indicating the polynomial order used in the discretization. Therefore, intuition is thus: the quadratic condition number term controlling the round-off error in a QR solve will be offset by a converging solution.

For instance, recall that $\mathrm{cond}(\widetilde{\mathsf{B}}) = \mathrm{cond}(\mathsf{A})^{1/2}$. Therefore, for many common FOS formulations, the condition number growth in $\widetilde{\mathsf{B}}$ is only $\mathcal{O}(h^{-1})$. Moreover, if the residual converges to zero as described above, then $\rho(\widetilde{\mathsf{B}}, \widetilde{\mathsf{l}}) \cdot \mathrm{cond}(\widetilde{\mathsf{B}})^2$ can be no worse than $\mathcal{O}(h^{-1})$. In such conventional scenarios, the numerical sensitivity of the least-squares solution is controlled only by the inverse of the mesh size and will be far more accurate than any normal equation approach!

Precisely, in the typical FOS scenario, we expect a QR-based algorithm will deliver an error bound of

$$\|\mathbf{e}\|_2 \leq \epsilon_{mach.} \|\mathbf{u}\|_2 C h^{-1} \,, \tag{3.15}$$

where $\mathbf{e}$ is the round-off error in the computation of the least-squares solution, $\epsilon_{mach.}$ is machine precision, and $C$ is a mesh-independent constant.

*An ill-conditioned Gram matrix.* Unfortunately, although the QR is guaranteed to deliver a more accurate solution than solving with the normal equation, there is potentially still a concealed obstacle. As many authors have pointed out, explicitly forming a product of two matrices before solving a least-squares problem posed with the matrix product $\widetilde{\mathsf{B}} = \mathsf{W}^{-1}\mathsf{B}$ is still not backwards stable [8]. Indeed, even when $\mathsf{B}$ is sparse and the matrix $\mathsf{W}$ is diagonal—but extremely ill-conditioned—this can be a potential issue [9, 47]. Because of this concern, several algorithms exist in the numerical linear algebra literature for this very class of problems [60, 2, 46, 72]. Nevertheless, we believe that such precautions are unwarranted in all but the most exceptional problems that can be expected with a DLS method.

Recall the critical assumptions in Section 1 on the structure of the Gram matrix. We implicitly assume that the conditioning of $\mathsf{W}$ should not be badly behaved as the problem size grows. Indeed, in the cases where the Gram matrix comes from a DPG method (i.e. is block-diagonal) or from some other technique (perhaps a preconditioner estimate), $\mathsf{G}^{-1}$ can usually be generated using local element or patch information. This motivates us to assume that some measure of its local rank structure should stay constant or be uniformly bounded (with respect to the mesh size, $h$) as the mesh is refined. If this is true, preconditioning the Gram matrix—which itself often acts like a preconditioner—with its diagonal entries,

$$\mathsf{G} \mapsto \mathsf{D}^{-1/2}\mathsf{G}\mathsf{D}^{-1/2} \,, \quad \mathsf{B} \mapsto \mathsf{D}^{1/2}\mathsf{B} \,, \quad \mathsf{l} \mapsto \mathsf{D}^{1/2}\mathsf{l} \,, \tag{3.16}$$

where $\mathsf{D} = \mathrm{diag}(\mathsf{G})$, before locally factoring into $\mathsf{W}\mathsf{W}^*$ and performing back-substitution, should lead to robust results.

This diagonal preconditioning procedure has been more than adequate in all of our experiments thus far. However, another possibility for improving the condition of the Gram matrix is the modified Lagrangian approach suggested in [42]. Here,

$$\mathsf{G} \mapsto \mathsf{G} + \mathsf{B}\mathsf{S}\mathsf{B}^* \,, \tag{3.17}$$

where $\mathsf{S}$ is a user-defined Hermitian matrix. Because the equivalent equilibrium systems

$$\begin{bmatrix} \mathsf{G} & \mathsf{B} \\ \mathsf{B}^* & 0 \end{bmatrix} \begin{bmatrix} \mathbf{p} \\ \mathbf{u} \end{bmatrix} = \begin{bmatrix} \mathsf{l} \\ 0 \end{bmatrix} \quad \Longleftrightarrow \quad \begin{bmatrix} \mathsf{G} + \mathsf{B}\mathsf{S}\mathsf{B}^* & \mathsf{B} \\ \mathsf{B}^* & 0 \end{bmatrix} \begin{bmatrix} \mathbf{p} \\ \mathbf{u} \end{bmatrix} = \begin{bmatrix} \mathsf{l} \\ 0 \end{bmatrix}$$

have the same solution, a well-chosen matrix $\mathsf{S}$ may significantly improve the condition number of the $(1, 1)$ block. It may also be useful to handle the case of a singular $\mathsf{G}$. We have experimented with this technique in our experiments, however, we have not made any certain conclusions for an effective $\mathsf{S}$. For the reader willing to apply this technique in the DPG setting, we note that $\mathsf{B}\mathsf{S}\mathsf{B}^*$ should probably hold the same block diagonal structure as $\mathsf{G}$.

If these procedures fail, a third alternative approach is obviously to avoid condensing the system. This is similar to some of the saddle-point finite element methods mentioned in Section 1. For completeness, we now summarize the highlights of this strategy.

*Generalized least-squares.* For the most badly behaved problems, [8, 44] suggest beginning with the generalized least-squares problem

$$\text{(GLS)} \qquad \min_{\mathbf{u}\in\mathbb{F}^N,\,\mathbf{r}\in\mathbb{F}^M} \|\mathbf{r}\|_2 \quad \text{subject to} \quad \mathsf{B}\mathbf{u} + \mathsf{W}\mathbf{r} = \mathsf{l}\,, \quad \text{where} \quad \mathsf{W}\mathsf{W}^* = \mathsf{G}\,.$$

For a direct method, the solution coefficients $\mathbf{u}$ are then suggested to be computed using a QR factorization approach which was first described in [58]. Although there are a couple of strict advantages of this idea—including that the Gram matrix $\mathsf{G}$ need no longer be invertible—the size of the resulting saddle-point system is much larger than the systems in the previously methods. Not to mention, the QR-based solution algorithm given in [58] is seemingly impractical for any reasonably large sparse system because it involves storing and applying large and probably dense orthogonal matrices. In [59], an efficient algorithm is proposed for the case that $\mathsf{G}$ is block diagonal, however, we are unaware of any multi-frontal implementations so have not explored it in our experiments.

Notably, this saddle-point approach is analogous to the finite element methods described in [20] and [25]. Moreover, although this generalized least-squares starting point may be too expensive to be practical for direct linear solvers, it may have benefits for iterative solution algorithms [74, 5, 43]. One such very promising technique in the DPG context, inspired by PDE-constrained optimization, is developed in [14] in a similar setting where $\mathsf{G}$ is not factored.

*Further discussion.* In choosing the best algorithm to solve the least-squares problem coming from a DLS method, many factors are important to consider. For instance, the normal equation have been demonstrated to be adequate when the methodology has been applied to many DPG problems [64, 65, 21, 67, 33, 37, 49, 62, 30, 40, 38, 35]. Indeed, in many reasonable circumstances, the round-off error in the solution from the associated linear solve cannot be expected to be nearly as large as the truncation error due to the finite element discretization. Direct solvers for Hermitian positive-definite systems are also usually several times faster than their associated QR-based counterparts for DLS problems.[1] Moreover, the normal equation, itself, generally requires less storage than the modified coefficient matrix and load vector.

Nevertheless, there are many anticipated circumstances where other least-squares methods would be especially useful. Archetypal examples include, but are not limited to, singular perturbation problems, problems with large material contrast, high-order PDEs, penalty methods, and nonlinear problems where the linear approximation may become singular or very ill-conditioned.

It is also important to mention that practical direct methods for constrained DLS methods have additional complications. Indeed, a naive direct QR-based algorithm for a problem

$$\text{(CDLS)} \qquad \min_{\mathbf{u}\in\mathbb{F}^N} \|\mathsf{B}\|_{\mathsf{G}^{-1}}\,, \quad \text{subject to} \quad \mathsf{C}\mathbf{u} = \mathsf{d}\,.$$

where both $\widetilde{\mathsf{B}}$ and $\mathsf{C}$ are sparse will probably not be practical. This is because such methods involve transforming the minimization problem above to an equivalent problem posed over the null space of $\mathsf{C}$. This sequence of operations will involve representing $\widetilde{\mathsf{B}}$ in a basis where it may no longer have a similar sparsity. The same issues appear for problems with inequality constraints [8, Section 6.8]. Therefore, for direct methods, (CDLS) requires either solving the condensed system

$$\begin{bmatrix} \mathsf{A} & \mathsf{C}^* \\ \mathsf{C} & 0 \end{bmatrix} \begin{bmatrix} \mathbf{u}^{\mathrm{opt}} \\ \mathbf{w}^{\mathrm{opt}} \end{bmatrix} = \begin{bmatrix} \mathsf{f} \\ \mathsf{d} \end{bmatrix},$$

or compromising on a penalized constraint and performing QR on (3.13).

[1]We make this statement based upon personal experience on a single-core machine with the available software at this time. Because most complexity arguments in the literature are usually pessimistic and proved for dense matrices (they do not incorporate the expected sparse matrix structure) and because the time-to-solution may depend heavily upon implementation or architecture, we believe this is more valid than predictions which could be made with well-known complexity results.

3.4. **Static condensation.** A common procedure which is often used to improve the solving time of linear systems is called static condensation. Here, the degrees of freedom associated with the element interior nodes (bubbles) are eliminated from the linear system. In practice, using a Schur complement procedure, small and independent blocks of the original stiffness matrix and load vectors are removed, and the original system is changed into a smaller-but-modified linear system with fewer unknowns. Often, this procedure of condensing, solving, and then recovering the global solution is much faster than solving the original system outright with standard means.

A similar procedure can be developed for overdetermined discrete least-squares problems without forming the normal equation—although it can be identified with static condensation therein (see Appendix B). To illustrate this application in practical DLS problems like (LS), consider (3.4) where we have separated $\mathsf{u} = [\mathsf{u}_{\text{bubb.}} \mid \mathsf{u}_{\text{interf.}}]^{\mathsf{T}}$ into bubble and interface components, respectively. Here, the bubble coefficients are assigned to functions which have support entirely within a single element, while the interface coefficients are assigned to basis functions with support across multiple elements or those which may only be defined at element interfaces. With this distinction, we may rewrite the minimization problem of (3.4) as a sequence of two independent minimizations:

$$\min_{\mathsf{u}\in\mathbb{F}^N} \mathcal{R}(\mathsf{u}) = \min_{\mathsf{u}_{\text{interf.}}\in\mathbb{F}^{N_{\text{interf.}}}} \min_{\mathsf{u}_{\text{bubb.}}\in\mathbb{F}^{N_{\text{bubb.}}}} \mathcal{R}(\mathsf{u}) = \min_{\mathsf{u}_{\text{bubb.}}\in\mathbb{F}^{N_{\text{bubb.}}}} \min_{\mathsf{u}_{\text{interf.}}\in\mathbb{F}^{N_{\text{interf.}}}} \mathcal{R}(\mathsf{u}) \,,$$

where we recall $\mathcal{R}(\mathsf{u}) = \big\|\widetilde{\mathsf{B}}\mathsf{u} - \widetilde{\mathsf{l}}\big\|_2^2$. Indeed, we may first solve for $\mathsf{u}^{\text{opt}}_{\text{interf.}}$ using an implicit expression for the unknown $\mathsf{u}^{\text{opt}}_{\text{bubb.}}$ in terms of the variable $\mathsf{u}_{\text{interf.}}$. Once $\mathsf{u}^{\text{opt}}_{\text{interf.}}$ has been computed, the second minimization problem,

$$\mathsf{u}^{\text{opt}}_{\text{bubb.}} = \operatorname*{arg\,min}_{\mathsf{u}_{\text{bubb.}}\in\mathbb{F}^{N_{\text{bubb.}}}} \mathcal{R}\big([\mathsf{u}_{\text{bubb.}} \mid \mathsf{u}^{\text{opt}}_{\text{interf.}}]\big) \,,$$

can be solved locally.

Representing the distinction in coefficients, $\mathsf{u} = [\mathsf{u}_{\text{bubb.}} \mid \mathsf{u}_{\text{interf.}}]^{\mathsf{T}}$, with the operator $\widetilde{\mathsf{B}}$ as well, we define two sub-matrices $\widetilde{\mathsf{B}}_{\text{bubb.}}$ and $\widetilde{\mathsf{B}}_{\text{interf.}}$ where $\widetilde{\mathsf{B}} = \big[\widetilde{\mathsf{B}}_{\text{bubb.}} \mid \widetilde{\mathsf{B}}_{\text{interf.}}\big]$. Naturally, the dimensions of $\widetilde{\mathsf{B}}_{\text{bubb.}}$ are $M \times N_{\text{bubb.}}$ and the dimensions of $\widetilde{\mathsf{B}}_{\text{interf.}}$ are $M \times N_{\text{interf.}}$. Recall that we are solving

$$\widetilde{\mathsf{B}}_{\text{bubb.}}\mathsf{u}_{\text{bubb.}} + \widetilde{\mathsf{B}}_{\text{interf.}}\mathsf{u}_{\text{interf.}} \overset{\text{LS}}{=} \widetilde{\mathsf{l}}\,, \tag{3.18}$$

where $\overset{\text{LS}}{=}$ indicates equality in the least-squares sense. If $\mathsf{u}_{\text{interf.}}$ were fixed, then the bubble contribution of the least-squares solution would be represented as

$$\mathsf{u}^{\text{opt}}_{\text{bubb.}} = \widetilde{\mathsf{B}}^{+}_{\text{bubb.}}\Big(\widetilde{\mathsf{l}} - \widetilde{\mathsf{B}}_{\text{interf.}}\mathsf{u}_{\text{interf.}}\Big) \,, \tag{3.19}$$

where $\widetilde{\mathsf{B}}^{+}_{\text{bubb.}} = \Big(\widetilde{\mathsf{B}}^{*}_{\text{bubb.}}\widetilde{\mathsf{B}}_{\text{bubb.}}\Big)^{-1}\widetilde{\mathsf{B}}^{*}_{\text{bubb.}}$ is the pseudoinverse of $\widetilde{\mathsf{B}}_{\text{bubb.}}$. Exploiting this expression for $\mathsf{u}^{\text{opt}}_{\text{bubb.}}$, we can rewrite (3.18) as a least-squares problem entirely in the variable $\mathsf{u}_{\text{interf.}}$ with $\mathsf{u}_{\text{bubb.}} = \mathsf{u}^{\text{opt}}_{\text{bubb.}}$ fixed, *viz.*,

$$(\mathrm{I} - \mathsf{P}_{\text{bubb.}})\widetilde{\mathsf{B}}_{\text{interf.}}\mathsf{u}_{\text{interf.}} \overset{\text{LS}}{=} (\mathrm{I} - \mathsf{P}_{\text{bubb.}})\widetilde{\mathsf{l}}\,, \tag{3.20}$$

where I is the identity matrix and $\mathsf{P}_{\text{bubb.}} = \widetilde{\mathsf{B}}_{\text{bubb.}}\Big(\widetilde{\mathsf{B}}^{*}_{\text{bubb.}}\widetilde{\mathsf{B}}_{\text{bubb.}}\Big)^{-1}\widetilde{\mathsf{B}}^{*}_{\text{bubb.}}$. The reader may observe that $\mathsf{P}^2_{\text{bubb.}} = \mathsf{P}_{\text{bubb.}}$ and likewise for $\mathsf{P}_{\text{interf.}} = \mathrm{I} - \mathsf{P}_{\text{bubb.}}$. These two matrices can be identified with orthogonal projectors onto the space of bubble coefficients and interface coefficients, respectively.

Applying $\mathsf{P}_{\text{interf.}}$ to $\widetilde{\mathsf{B}}_{\text{interf.}}$ can be performed in parallel because $\widetilde{\mathsf{B}}_{\text{bubb.}}$ will always be block diagonal. Then, without actually computing $\mathsf{u}^{\text{opt}}_{\text{bubb.}}$, expression (3.20) can be used to solve for the optimal interface coefficients,

$$\mathsf{u}^{\text{opt}}_{\text{interf.}} = \operatorname*{arg\,min}_{\mathsf{u}_{\text{interf.}}\in\mathbb{F}^{N_{\text{interf.}}}} \mathcal{R}\big([\mathsf{u}^{\text{opt}}_{\text{bubb.}} \mid \mathsf{u}_{\text{interf.}}]\big) \,.$$

Afterwards, once $\mathsf{u}^{\text{opt}}_{\text{interf.}}$ are available, (3.19) can be used to extract $\mathsf{u}^{\text{opt}}_{\text{bubb.}}$, locally.

Before closing, we remark that with a QR-decomposition of $\widetilde{\mathsf{B}}_{\text{bubb.}}$,

$$\widetilde{\mathsf{B}}_{\text{bubb.}} = \Big[\mathsf{Q}_{\text{bubb.}} \mid \mathsf{Q}_{\text{interf.}}\Big]\begin{bmatrix}\mathsf{R}_{\text{bubb.}}\\ 0\end{bmatrix}$$

where $\left[\mathsf{Q}_{\text{bubb.}} \mid \mathsf{Q}_{\text{interf.}}\right]$ is unitary and $\mathsf{R}_{\text{bubb.}}$ is upper-triangular, (3.19) can be rewritten as

$$\mathbf{u}^{\text{opt}}_{\text{bubb.}} = \mathsf{R}^{-1}_{\text{bubb.}}\mathsf{Q}^{*}_{\text{bubb.}}\left(\widetilde{\mathsf{l}} - \widetilde{\mathsf{B}}_{\text{interf.}}\mathbf{u}_{\text{interf.}}\right).$$

Meanwhile, $\mathsf{P}_{\text{bubb.}} = \mathsf{Q}_{\text{bubb.}}\mathsf{Q}^{*}_{\text{bubb.}}$ and $\mathsf{P}_{\text{interf.}} = \mathrm{I} - \mathsf{Q}_{\text{bubb.}}\mathsf{Q}^{*}_{\text{bubb.}} = \mathsf{Q}_{\text{interf.}}\mathsf{Q}^{*}_{\text{interf.}}$.

Exploiting the QR-decomposition will minimize the round-off error introduced in modifying the original system as well as in recovering $\mathbf{u}^{\text{opt}}_{\text{bubb.}}$ and we suggest static condensation routines for DLS problems of the form (LS) be implemented using it.

## 4. Assembly

In this section, we describe the construction of the global linear systems for DPG methods. As it was outlined in Section 3.3, we are primarily interested in two different procedures: (1) directly assembling the normal equation; and (2) assembling the overdetermined system. To date, forming the normal equation has been the primary assembly procedure for DPG linear systems. The main advantages of this approach is that the assembly algorithm is identical to all traditional conforming finite element methods and it will involve the least storage. Moreover, many efficient direct and iterative solvers specialized for Hermitian/symmetric positive definite systems can be employed. Nonetheless, there is an important disadvantage: the condition number of this global stiffness matrix $\mathsf{A}$ is the square of the condition number of the alternative, which is the (global) enriched stiffness matrix $\widetilde{\mathsf{B}}$. The enriched stiffness matrix is, however, rectangular, and so other solvers—generally more expensive solvers—have to be used to solve the overdetermined linear system it is involved with. Be that as it may, this second approach can be applied to ill-conditioned problems, where forming the normal equation becomes an unsatisfactory option.

We proceed by giving a brief description of the two different assembly procedures.

### 4.1. The normal equation.

The assembly of the normal equation for the DPG method can easily be incorporated into any finite element code supporting exact sequence conforming shape functions [3, 36]. Recall that the DPG Gram matrix $\mathsf{G}$ is block diagonal. Therefore, it can be inverted element-wise and, therefore, $\mathsf{G}^{-1}$ is also block diagonal. Let $\mathsf{B}_K$ denote the enriched stiffness matrix for element $K$ and $\mathsf{l}_K$ the corresponding load vector. Additionally, let $\mathsf{G}_K$ be the element Gram matrix. Then, the DPG element stiffness matrix $\mathsf{A}_K$ and load vector $\mathsf{f}_K$ are given by $\mathsf{A}_K = \mathsf{B}^*_K\mathsf{G}^{-1}_K\mathsf{B}_K$ and $\mathsf{f}_K = \mathsf{B}^*_K\mathsf{G}^{-1}_K\mathsf{l}_K$, respectively. Using the Cholesky factorization of $\mathsf{G}_K = \mathsf{L}_K\mathsf{L}^*_K$, we have:

$$\mathsf{A}_K = \mathsf{B}^*_K\mathsf{G}^{-1}_K\mathsf{B}_K = (\mathsf{L}^{-1}_K\mathsf{B}_K)^*(\mathsf{L}^{-1}_K\mathsf{B}_K) \quad \text{and} \quad \mathsf{f}_K = \mathsf{B}^*_K\mathsf{G}^{-1}_K\mathsf{l}_K = (\mathsf{L}^{-1}_K\mathsf{B}_K)^*(\mathsf{L}^{-1}_K\mathsf{l}_K).$$

We note that one may wish to precondition before the above operations, as in (3.16).

The computation of the element DPG stiffness matrix and load vector is given by Algorithm 1. The assembly

**Algorithm 1** Element stiffness matrix and load vector for the DPG normal equation

1: $\mathsf{L}_K \leftarrow$ Cholesky($\mathsf{G}_K$)
2: $\widetilde{\mathsf{B}}_K \leftarrow$ Triangular solve($\mathsf{L}_K\widetilde{\mathsf{B}}_K = \mathsf{B}_K$)
3: $\widetilde{\mathsf{l}}_K \leftarrow$ Triangular solve($\mathsf{L}_K\widetilde{\mathsf{l}}_K = \mathsf{l}$)
4: $\mathsf{A}_K \leftarrow \widetilde{\mathsf{B}}^*_K\widetilde{\mathsf{B}}_K$ // DPG element stiffness matrix
5: $\mathsf{f}_K \leftarrow \widetilde{\mathsf{B}}^*_K\widetilde{\mathsf{l}}_K$. // DPG element load vector of the global DPG stiffness matrix and load vector can be implemented by following the common algorithm of any standard finite element code [26, 31]. Note that there are two modifications that one should make to the element stiffness matrices before the global assembly: account for Dirichlet boundary conditions; and (optional) accommodate degrees of freedom associated to constrained nodes for adaptive mesh refinement (hanging nodes are possibly created after adaptive $h$-refinements). We refer the reader to [26, 31] for detailed discussion on both of these modifications.

After pre-processing the element stiffness matrices, $\mathsf{A}_K \mapsto \mathsf{A}_K^{\mathtt{mod}}$, one should proceed with static condensation, $\mathsf{A}_K^{\mathtt{mod}} \mapsto \mathsf{A}_K^{\mathtt{c}}$, to reduce the complexity of the global system. For these square and symmetric matrices, this operation is described in Appendix B.

Additionally, as in standard FEM, the assembly is driven by the so-called "local-to-global connectivity maps". These maps assign to the local degrees of freedom their corresponding global degrees of freedom. The construction of these maps is based on the "donor strategy" and is implemented as in [26].

A description of the assembly procedure is given in Algorithm 2.

**Algorithm 2** Assembly of DPG normal equation

Initialize global stiffness matrix and load vector $\mathsf{A}$ and $\mathsf{f}$.
**for** $K \leftarrow 1$ **to** $N_K$ **do** // for each element in the mesh
    Compute $\mathsf{A}_K$ and $\mathsf{f}_K$ // element stiffness matrix and load vector
    Compute $\mathsf{A}_K^{\mathtt{mod}}$ and $\mathsf{f}_K^{\mathtt{mod}}$ // modified element matrix and load vector
    Compute $\mathsf{A}_K^{\mathtt{c}}$ and $\mathsf{f}_K^{\mathtt{c}}$ // condensed element matrix and load vector
    Get $\mathtt{Con}_K$ // local-to-global connectivity map
    **for** $k_1 \leftarrow 1$ **to** $\mathtt{ndof}_K$ **do** // for each element degree of freedom (DOF)
        $i \leftarrow \mathtt{Con}_K(k_1)$ // global index for local DOF
        $\mathsf{f}(i) \leftarrow \mathsf{f}(i) + \mathsf{f}_K(k_1)$ // accumulate for the global load vector
        **for** $k_2 \leftarrow 1$ **to** $\mathtt{ndof}_K$ **do** // for each element DOF
            $j \leftarrow \mathtt{Con}_K(k_2)$ // global index for local DOF
            $\mathsf{A}(i,j) \leftarrow \mathsf{A}(i,j) + \mathsf{A}_K(k_1,k_2)$ // accumulate for the global stiffness matrix

4.2. **The overdetermined system.** Constructing the global overdetermined system requires some modifications to the assembly algorithms above. First, in order to deliver *rectangular* element stiffness matrices $\widetilde{\mathsf{B}}$ and load vectors $\widetilde{\mathsf{l}}$, one should only perform the first three steps of Algorithm 1. Note that the column size of the element stiffness matrix $\widetilde{\mathsf{B}}$ corresponds to the number of trial degrees of freedom and the row size to the number of test degrees of freedom. Similarly, the size of the load vector $\widetilde{\mathsf{l}}$ corresponds to the number of the test degrees of freedom.

As with square stiffness matrices, after the rectangular element matrices and load vectors have been computed, they need to be modified in order to accommodate Dirichlet boundary conditions and constrained nodes. The Dirichlet boundary conditions can be accounted for, $\widetilde{\mathsf{B}} \mapsto \widetilde{\mathsf{B}}_{\text{hom.}}$ and $\widetilde{\mathsf{l}} \mapsto \widetilde{\mathsf{l}}_{\text{lift}}$, through the operation described in Section 3.2. This, however, is performed locally, like in the assembly algorithm for the normal equation, so that $\widetilde{\mathsf{B}}_K \mapsto \widetilde{\mathsf{B}}_K^{\mathtt{mod}}$ and $\widetilde{\mathsf{l}}_K \mapsto \widetilde{\mathsf{l}}_K^{\mathtt{mod}}$.

For constrained nodes, the procedure is similar to the one for the normal equation, with the difference being that now the modifications are performed only on the trial space because the test space is broken. Note that there are no modifications needed for the load vector.

The final local step is static condensation, $\widetilde{\mathsf{B}}_K^{\mathtt{mod}} \mapsto \widetilde{\mathsf{B}}_K^{\mathtt{c}}$ and $\widetilde{\mathsf{l}}_K^{\mathtt{mod}} \mapsto \widetilde{\mathsf{l}}_K^{\mathtt{c}}$ (see Section 3.4).

The global assembly algorithm then proceeds in a similar manner as for the normal equation. However, there is one important difference because the test space is broken: the need for accumulation of the contributions from different elements in both the (global) enriched stiffness matrix $\widetilde{\mathsf{B}}$ and the (global) enriched load vector $\widetilde{\mathsf{l}}$ has been eliminated. Therefore, in the global stiffness matrix, every row is independent. Note that this allows for a fully parallel assembly algorithm.

This entire global assembly procedure is summarized in Algorithm 3.

Note that, in practice, both Algorithms 2 and 3 are modified to exploit sparsity of the stiffness matrices, $\mathsf{A}$ and $\widetilde{\mathsf{B}}$, depending on what solution package is used for the linear solve.

4.3. **Comparison.** For the one-dimensional DPG ultraweak formulation of Poisson equation (see Section 5), Figures 4.1 and 4.2 depict the global DPG stiffness matrix for the normal equation and the overdetermined system, respectively. The mesh used consisted of ten quadratic elements and the order of approximation for the

**Algorithm 3** Assembly of DPG overdetermined system

Initialize global stiffness matrix and load vector $\widetilde{\mathsf{B}}$ and $\widetilde{\mathsf{f}}$.
Initialize global test DOF counter $i$.
**for** $K \leftarrow 1$ **to** $N_K$ **do** // for each element in the mesh
  Compute $\widetilde{\mathsf{B}}_K$ and $\widetilde{\mathsf{l}}_K$ // element stiffness matrix and load vector
  Compute $\widetilde{\mathsf{B}}_K^{\mathrm{mod}}$ and $\widetilde{\mathsf{l}}_K^{\mathrm{mod}}$ // modified element matrix and load vector
  Compute $\widetilde{\mathsf{B}}_K^{\mathrm{c}}$ and $\widetilde{\mathsf{l}}_K^{\mathrm{c}}$ // condensed element matrix and load vector
  Get $\mathtt{Con}_K$ // local-to-global connectivity map
  **for** $k_1 \leftarrow 1$ **to** $\mathtt{ndofT}_K$ **do** // for each element test degree of freedom
    $i \leftarrow i+1$ // global test DOF counter
    $\widetilde{\mathsf{l}}(i) \leftarrow \widetilde{\mathsf{l}}_K(k_1)$ // fill in the global load vector
    **for** $k_2 \leftarrow 1$ **to** $\mathtt{ndof}_K$ **do** // for each element DOF
      $j \leftarrow \mathtt{Con}_K(k_2)$ // global index for local DOF
      $\widetilde{\mathsf{B}}(i,j) \leftarrow \widetilde{\mathsf{B}}_K(k_1,k_2)$ // fill in the global stiffness matrix

enriched test space was three. In a broken ultraweak formulation, continuity is enforced with the introduction of new interface unknowns (Lagrange multipliers), usually called fluxes and traces, that are traces of dual energy spaces on the mesh skeleton. This explains the structure of the matrix for the total system seen in Figure 4.1 (A). After static condensation of the interior degrees of freedom, the resulting linear system involved only the interface unknowns, and, therefore, the matrix in Figure 4.1 (B) consists of only one band of overlapping blocks.

The situation is slightly different in the case of the overdetermined system (see Figure 4.2 (A)). For instance, because the test space is broken, there is no overlap between rows. Similar to the case of the normal equation, static condensation led to a linear system involving only the interfaces unknowns. However, here, the size of test space remained the same and therefore only the column dimension was reduced (see Figure 4.2 (B)).

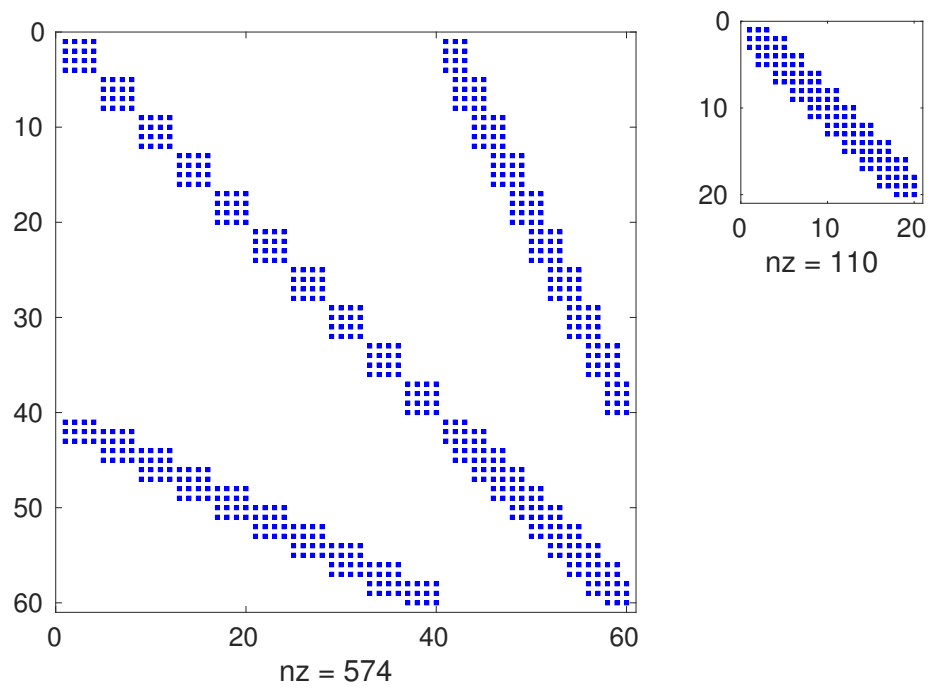

(A) Total system (B) Condensed system

FIGURE 4.1. DPG stiffness matrix A.

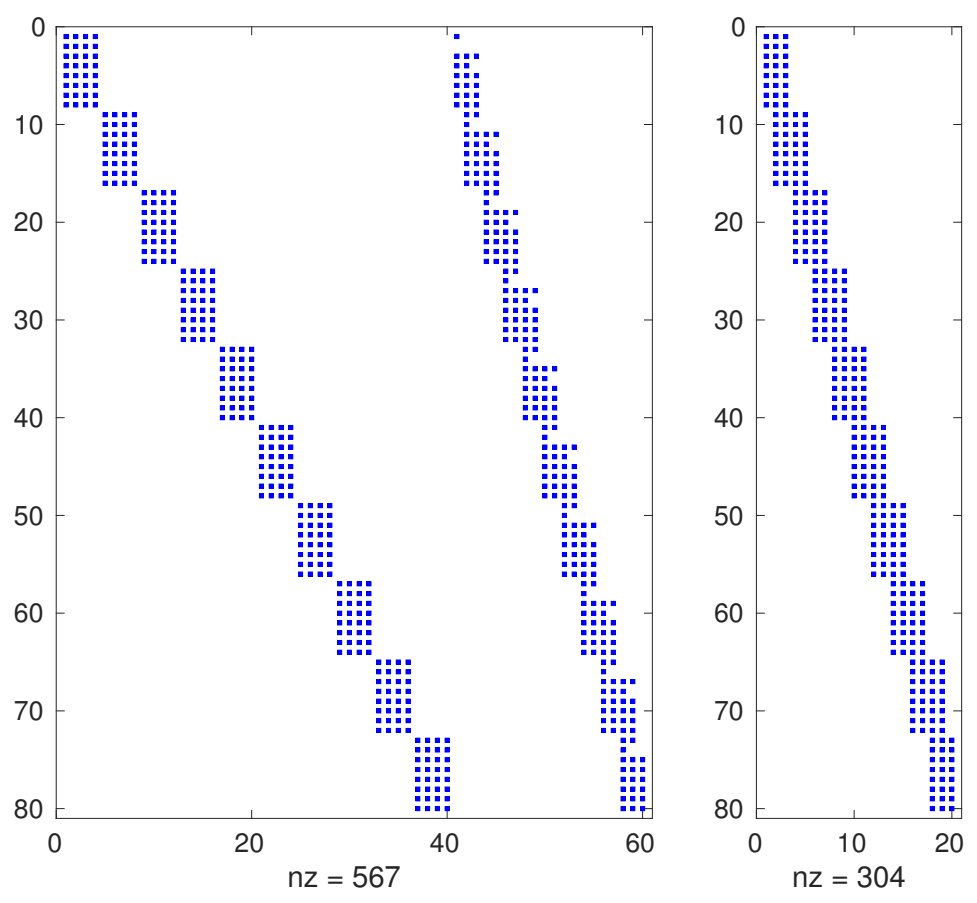

(A) Total system (B) Condensed system

FIGURE 4.2. DPG stiffness matrix $\widetilde{\mathsf{B}}$.

## 5. CONDITIONING

5.1. **Robustness.** Observe that if $\mathcal{B}(\mathcal{U}_h) \subset \mathcal{R}_{\mathcal{V}}\, \mathcal{V}_r$ then $\mathcal{B}(\mathcal{U}_h) \perp \mathcal{R}_{\mathcal{V}}\, \mathcal{V}_r^{\perp}$ for $\mathcal{V} = \mathcal{V}_r \oplus \mathcal{V}_r^{\perp}$. Therefore,

$$\mathcal{R}_{\mathcal{V}}^{-1}\, \mathcal{B}(\mathcal{U}_h) = \big(\mathcal{R}_{\mathcal{V}_r}^{-1} \oplus \mathcal{R}_{\mathcal{V}_r^{\perp}}^{-1}\big)\mathcal{B}(\mathcal{U}_h) = \mathcal{R}_{\mathcal{V}_r}^{-1}\, \mathcal{B}(\mathcal{U}_h) \oplus 0 = \mathcal{R}_{\mathcal{V}_r}^{-1}\, \mathcal{B}(\mathcal{U}_h)\,.$$

This indicates that, in special circumstances, no error is introduced from discretizing the Riesz map $\mathcal{R}_\mathcal{V}$. Using this result, rewriting (2.6), and proceeding similarly with (3.1), yields

$$\left\langle \mathcal{B}\mathfrak{u}_{h,r}^{\mathrm{opt}} - \ell, \mathcal{R}_{\mathcal{V}_r}^{-1} \mathcal{B}\delta\mathfrak{u}_h \right\rangle = 0 = \left\langle \mathcal{B}\mathfrak{u}_h^{\mathrm{opt}} - \ell, \mathcal{R}_{\mathcal{V}}^{-1} \mathcal{B}\delta\mathfrak{u}_h \right\rangle = \left\langle \mathcal{B}\mathfrak{u}_h^{\mathrm{opt}} - \ell, \mathcal{R}_{\mathcal{V}_r}^{-1} \mathcal{B}\delta\mathfrak{u}_h \right\rangle . \tag{5.1}$$

Thus, if $\mathcal{B}(\mathcal{U}_h) \subset \mathcal{R}_\mathcal{V}\, \mathcal{V}_r$, then $\mathfrak{u}_h^{\mathrm{opt}} = \mathfrak{u}_{h,r}^{\mathrm{opt}}$; meaning that the discrete solution is actually the same solution which minimizes the continuous residual. In the particular case of least-squares methods, as in (2.8), this implies

$$(\mathsf{A}_{\mathrm{LS}})_{ij} = (\mathcal{L}\mathfrak{u}_j, \mathcal{L}\mathfrak{u}_i)_{L^2(\Omega)} = \left\langle \mathcal{B}'\, \mathcal{R}_{L^2(\Omega)}^{-1}\, \mathcal{B}\mathfrak{u}_j, \mathfrak{u}_i \right\rangle = \left\langle \mathcal{B}'\, \mathcal{R}_{\mathcal{V}_r}^{-1}\, \mathcal{B}\mathfrak{u}_j, \mathfrak{u}_i \right\rangle = \mathsf{A}_{ij} . \tag{5.2}$$

Hence, even the monolithic coefficient matrix of a least-squares method, $\mathsf{A}_{\mathrm{LS}}$, is exactly the same as that of a DLS discretization in the normal equation setting, provided $\mathcal{B}(\mathcal{U}_h) \subset \mathcal{R}_\mathcal{V}\, \mathcal{V}_r$. However, only with the DLS discretization is the factorization of the matrix $\mathsf{A} = \widetilde{\mathsf{B}}^*\widetilde{\mathsf{B}}$ readily available. One advantage being that with a QR algorithm, as discussed in Section 3.3, one could solve for exactly the same solution, but with less numerical sensitivity.

We will now compare the condition number behavior, numerical sensitivity, discretization error, and round-off error incurred in DPG methods. As in [63], in each of our experiments, the stiffness matrix was diagonally preconditioned:

$$\mathsf{A} \mapsto \mathsf{D}^{-1/2}\mathsf{A}\mathsf{D}^{-1/2} , \quad \mathsf{f} \mapsto \mathsf{D}^{-1/2}\mathsf{f} ,$$

or, equivalently,

$$\widetilde{\mathsf{B}} \mapsto \widetilde{\mathsf{B}}\mathsf{D}^{-1/2} , \quad \widetilde{\mathsf{l}} \mapsto \widetilde{\mathsf{l}} ,$$

where $\mathsf{D} = \mathrm{diag}(\mathsf{A})$. The cost of this procedure is computationally negligible and it is common practice to scale the matrix in this way before iterative solution methods. Meanwhile, it is performed implicitly in most direct solvers. Therefore, we presume no offense in this action. For additional perspective on several topics we do not cover, related to the condition number of DPG stiffness matrices, but with a focus on Stokes equation, we refer the interested reader to [63, Chapter 9].

*First-order system least squares.* For illustration, consider Poisson's equation in $\mathbb{R}^2$ with body force $f$ and Dirichlet boundary condition $u|_{\partial\Omega} = \hat{f}$. Here, the aim is to find $(u, \boldsymbol{\sigma}) \in (u^{\text{lift}} + H_0^1(\Omega)) \times \boldsymbol{H}(\mathrm{div}, \Omega)$ such that

$$\begin{aligned} -\operatorname{div} \boldsymbol{\sigma} + \alpha u &= f , \\ \boldsymbol{\sigma} - \operatorname{grad} u &= 0 , \end{aligned} \tag{5.3}$$

where, $\alpha = 0$, $f \in L^2$, and $u^{\text{lift}} \in H^1(\Omega)$ is an extension to $\Omega$ where $\mathrm{tr}_{H^1(\Omega)} u^{\text{lift}} = \hat{f} \in H^{1/2}(\partial\Omega)$. These equations can be solved using the first-order system least-squares (FOSLS) method [16, 17] by directly discretizing the trial space $\mathcal{U}^{\text{hom.}} = H_0^1(\Omega) \times \boldsymbol{H}(\mathrm{div}, \Omega)$, as in (2.8). The appropriate discretization $\mathcal{U}_h^{\text{hom.}} \subset \mathcal{U}^{\text{hom.}}$ should be locally derived from an exact sequence of spaces of polynomial order $p$, such as the ones found in [36] for triangles and quadrilaterals. Indeed, for a master square, $(0, 1)^2$, consider the spaces $W^p = \mathcal{Q}^{p,p}$, $\boldsymbol{V}^p = \mathcal{Q}^{p,p-1} \times \mathcal{Q}^{p-1,p}$ and $Y^p = \mathcal{Q}^{p-1,p-1}$, where $\mathcal{Q}^{p,q} = \mathcal{P}^p(x) \otimes \mathcal{P}^q(y)$. Thus, for a quadrilateral mesh, a compatible basis $\mathfrak{U}_h^{\text{hom.}}$ of $\mathcal{U}_h^{\text{hom.}} \subset \mathcal{U}^{\text{hom.}}$ can be constructed, locally, at each element $K \in \mathcal{T}$, as the pullback of $W^p \times \boldsymbol{V}^p$ from the master square to $K$.

On the other hand, solving (5.3) using a DLS discretization requires defining a bilinear form with a test space $\mathcal{V}$ and its corresponding discretization $\mathcal{V}_r \subset \mathcal{V}$. Thus, dropping the superscript-hom. from now on, the goal is that of finding $(u, \boldsymbol{\sigma}) \in H_0^1(\Omega) \times \boldsymbol{H}(\mathrm{div}, \Omega) = \mathcal{U}$ such that

$$b\big((u, \boldsymbol{\sigma}), (v, \boldsymbol{\tau})\big) = \ell\big((v, \boldsymbol{\tau})\big) , \quad \text{for all } (v, \boldsymbol{\tau}) \in L^2(\Omega) \times \boldsymbol{L}^2(\Omega) = \mathcal{V} ,$$

where

$$\begin{aligned} b\big((u, \boldsymbol{\sigma}), (v, \boldsymbol{\tau})\big) &= -(\operatorname{div} \boldsymbol{\sigma}, v)_\Omega + (\alpha u, v)_\Omega + (\boldsymbol{\sigma}, \boldsymbol{\tau})_\Omega - (\operatorname{grad} u, \boldsymbol{\tau})_\Omega , \\ \ell\big((v, \boldsymbol{\tau})\big) &= (f, v)_\Omega - (\alpha u^{\text{lift}}, v)_\Omega + (\operatorname{grad} u^{\text{lift}}, \boldsymbol{\tau})_\Omega . \end{aligned} \tag{5.4}$$

And, where $(\cdot,\cdot)_\Omega = (\cdot,\cdot)_{L^2(\Omega)}$, $\boldsymbol{L}^2(\Omega) = L^2(\Omega) \times L^2(\Omega)$ (with the standard norm $\langle \mathcal{R}_\mathcal{V}\, \mathfrak{v}, \mathfrak{v} \rangle = \|\mathfrak{v}\|_\mathcal{V}^2 = \|\mathfrak{v}\|_{\boldsymbol{L}^2(\Omega)}^2$ for all $\mathfrak{v} = (v, \boldsymbol{\tau}) \in \mathcal{V}$), and $(u + u^{\text{lift}}, \boldsymbol{\sigma})$ is the solution to (5.3). This is known as the (first-order)

"strong variational formulation" of the Poisson equation [19]. A DLS discretization can be found by using the same finite-dimensional spaces $\mathcal{U}_h$ as in the FOSLS method above, as well as an $L^2$-conforming test space $\mathcal{V}_r$ derived from an exact sequence of order $p+\Delta p$. Here, locally, $\mathcal{V}_r$ would be a pullback of $(Y^{p+\Delta p})^3$. Notice that, at the master element level, $\operatorname{div}\boldsymbol{\sigma} \in Y^p$, while $\boldsymbol{\sigma} - \operatorname{grad} u \in (Y^{p+1})^2$. Thus, it follows that $\mathcal{B}(\mathcal{U}_h) \subset \mathcal{R}_\mathcal{V}\,\mathcal{V}_r$ whenever $\Delta p \geq 1$.

Let $\Omega = (0,1)^2$ be partitioned by a uniform quadrilateral mesh of side length $h$, and assume $u(x,y) = \sin(\pi x)\sin(\pi y)$. Define $\hat{f} = 0$ and $f = -\operatorname{div}(\operatorname{grad} u)$. With the lift $u^{\text{lift}} = 0$, we separately solved (5.3) using both least-squares (i.e. FOSLS) and DLS methods (with $\Delta p \geq 1$). Moreover, because of the uniformity of the mesh and from using $L^2$-orthogonal shape functions at the master element level [36], we produced a diagonal Gram matrix, $\mathsf{G}_{ij} = (\mathfrak{v}_i, \mathfrak{v}_j)_{\boldsymbol{L}^2(\Omega)}$.

The condition number we computed for each of the three possible stiffness matrices is presented in Figure 5.1. Note that the condition number of the FOSLS stiffness matrix $\mathsf{A}_{\text{LS}}$ was verified to grow as $\mathcal{O}(h^{-2})$ (same as with a DLS stiffness matrix $\mathsf{A}$), but the condition number of the DLS enriched stiffness matrix $\widetilde{\mathsf{B}}$ was verified to grow only as $\mathcal{O}(h^{-1})$.

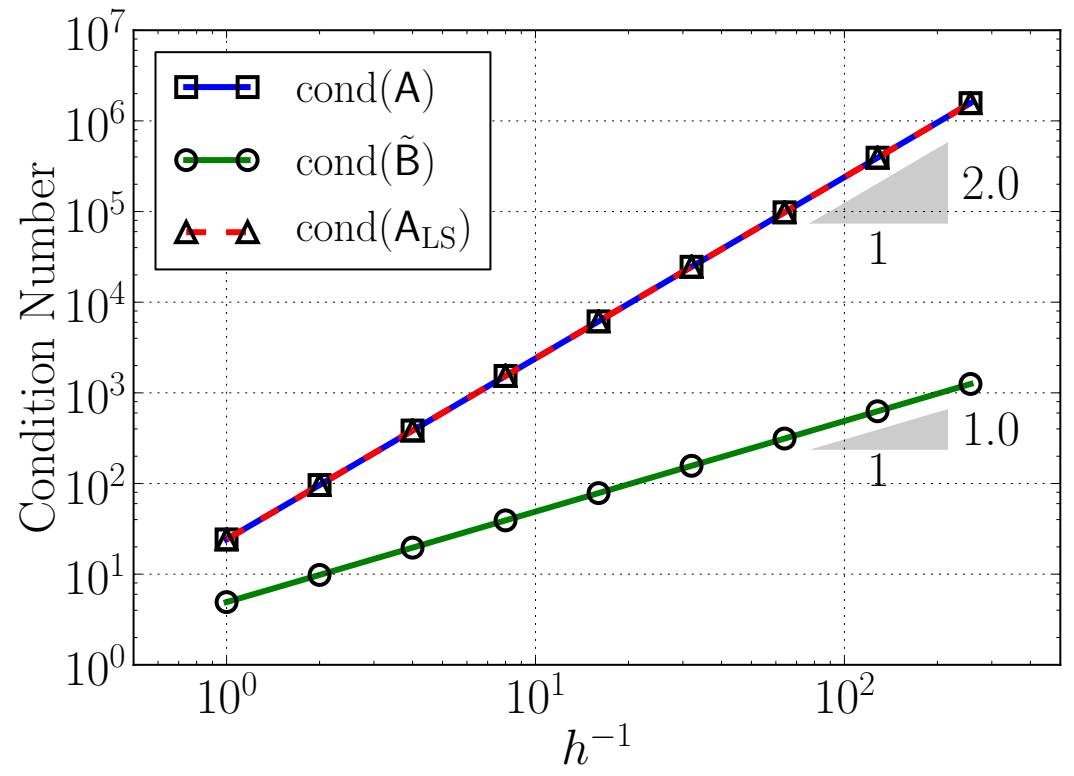

FIGURE 5.1. Comparison with Poisson's equation of the condition number of the FOSLS stiffness matrix $\mathsf{A}_{\text{LS}}$ to the condition number of the DLS stiffness matrices $\mathsf{A}$ and $\widetilde{\mathsf{B}}$ coming from the strong formulation, (5.4). Exact sequence polynomial order $p = 2$ for $\mathcal{U}_h$ and, in the DLS setting, $\Delta p = 1$. Observe that $\operatorname{cond}(\mathsf{A}_{\text{LS}}) = \operatorname{cond}(\mathsf{A}) = \operatorname{cond}(\widetilde{\mathsf{B}})^2$. All reported results are for the statically condensed and diagonally preconditioned matrices.

Recall (2.4) and (3.1). The FOSLS solution would be $\mathfrak{u}_h^{\text{opt}}$, while that coming from the DLS discretization would be $\mathfrak{u}_{h,r}^{\text{opt}}$. As expected from (5.2), numerical results also confirmed that when $\Delta p \geq 1$, $\mathfrak{u}_h^{\text{opt}} = \mathfrak{u}_{h,r}^{\text{opt}}$ and $\mathsf{A}_{\text{LS}} = \mathsf{A}$, up to floating-point precision.

When $\mathcal{B}(\mathcal{U}_h) \not\subset \mathcal{R}_\mathcal{V}\,\mathcal{V}_r$, the solutions $\mathfrak{u}_h^{\text{opt}}$ and $\mathfrak{u}_{h,r}^{\text{opt}}$, and matrices $\mathsf{A}_{\text{LS}}$ and $\mathsf{A}$ will no longer be equal. Naturally, however, the distance between $\mathfrak{u}_{h,r}^{\text{opt}}$ to $\mathfrak{u}_h^{\text{opt}}$ is expected to decrease as $\mathcal{V}_r$ is enriched (i.e. as $\Delta p$ is increased). Specifically, if the enriched test spaces are nested, $\mathcal{V}_{r_1} \subsetneq \mathcal{V}_{r_2} \subsetneq \mathcal{V}_{r_3} \subsetneq \cdots$, we expect $\|\mathfrak{u}_h^{\text{opt}} - \mathfrak{u}_{h,r_k}^{\text{opt}}\|_\mathcal{U} \to 0$ as $k \in \mathbb{N}$ increases. Indeed, this was observed when considering $\alpha(x,y) = \sin(\pi x)\sin(\pi y)$ in (5.3) and (5.4) and comparing the FOSLS solution, $\mathfrak{u}_h^{\text{opt}}$, to the DLS normal equation solution, $\mathfrak{u}_{h,\Delta p}^{\text{opt}}$, for increasing values of $\Delta p$. The results in Figure 5.2 (A) show that the rate of $h$-convergence between the two discrete solutions grows with $\Delta p$. Moreover, Figure 5.2 (B) shows that the matrix $\mathsf{A}$ converges to $\mathsf{A}_{\text{LS}}$. These numerical results suggest that the error incurred in discretizing $\mathcal{V}$ to $\mathcal{V}_r$ can be made very small, and we expect this to be true even with non-trivial variational formulations (i.e. after integrating by parts).

*Primal DPG and Bubnov-Galerkin methods.* In order to explore less trivial variational formulations, consider the broken primal formulation of the Poisson equation in (5.3) (with $\alpha = 0$) [28]. In this setting, we seek a solution $(u, \hat{\sigma}_n) \in H_0^1(\Omega) \times H^{-1/2}(\partial\mathcal{T}) = \mathcal{U}$ (see [19, 48] for definitions of mesh-trace Sobolev spaces) such that

$$b\big((u, \hat{\sigma}_n), v\big) = \ell(v)\,, \quad \text{for all } v \in H^1(\mathcal{T}) = \mathcal{V}\,,$$

where

$$\begin{aligned} b\big((u, \hat{\sigma}_n), v\big) &= (\operatorname{grad} u, \operatorname{grad} v)_{\mathcal{T}} - \langle \hat{\sigma}_n, v \rangle_{\partial\mathcal{T}}\,, \\ \ell(v) &= (f, v)_{\mathcal{T}} - (\operatorname{grad} u^{\text{lift}}, \operatorname{grad} v)_{\mathcal{T}}\,. \end{aligned} \tag{5.5}$$

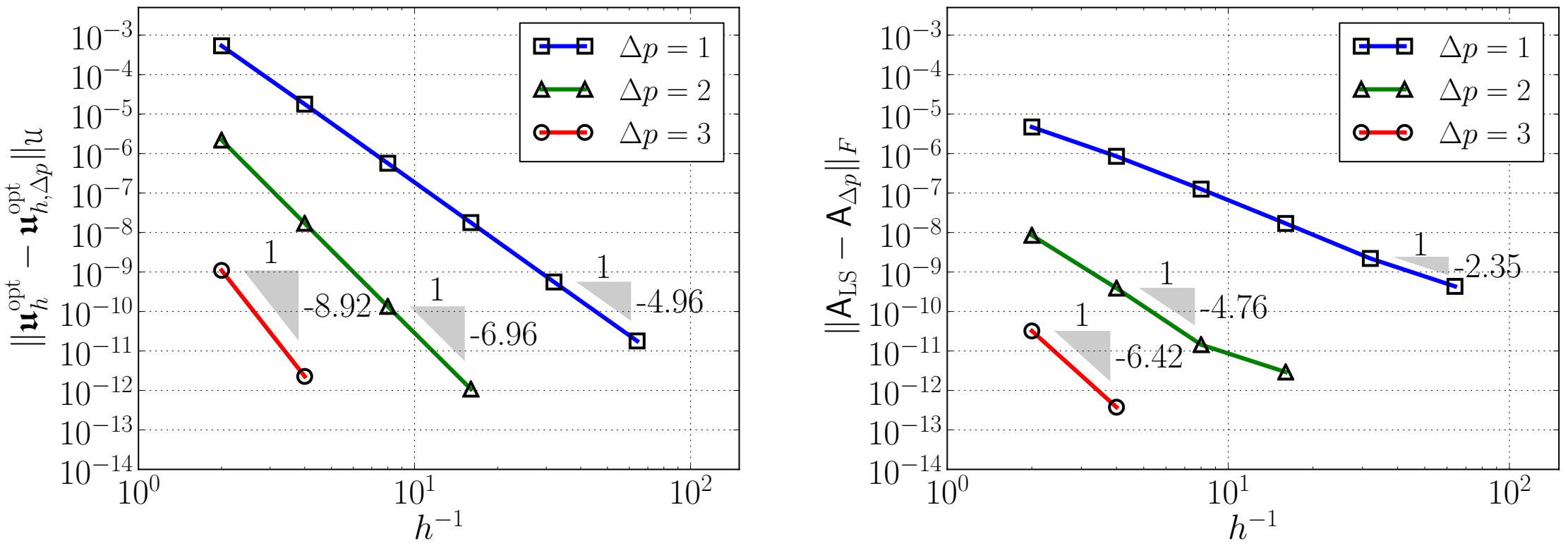

(A) Solution difference

(B) Matrix difference

FIGURE 5.2. Distance between the FOSLS and DLS solutions, $\mathbf{u}_h^{\text{opt}}$ and $\mathbf{u}_{h,\Delta p}^{\text{opt}}$, and the stiffness matrices $\mathsf{A}_{\text{LS}}$ and $\mathsf{A}$. The solution $\mathbf{u}_{h,\Delta p}^{\text{opt}}$ was computed, in each experiment, by way of the normal equation (NE). Exact sequence polynomial order $p = 2$ was used for $\mathcal{U}_h$, for all computations. Note that $\|(u, \boldsymbol{\sigma})\|_{\mathcal{U}}^2 = \|u\|_{H^1}^2 + \|\boldsymbol{\sigma}\|_{H(\text{div})}^2$ and $\|\cdot\|_F^2$ is the Frobenius norm.

Here, $H^1(\mathcal{T})$ is a *broken* Sobolev space with norm $\|v\|_{H^1(\mathcal{T})}^2$ [19], and the restriction of each member of this space to any single $K \in \mathcal{T}$ is in $H^1(K)$. Likewise, $(\cdot,\cdot)_{\mathcal{T}} = \sum_{K\in\mathcal{T}} (\cdot,\cdot)_K$ and similarly with $\langle\cdot,\cdot\rangle_{\partial\mathcal{T}}$. This second pairing, $\langle\cdot,\cdot\rangle_{\partial\mathcal{T}}$, can be understood, intuitively, as a mesh-boundary integral, however, the inquisitive reader may wish to examine [19, 48] for further detail. The lift $u^{\text{lift}} \in H^1(\Omega)$ in (5.5) is identical to that from the least-squares setting presented before.

For discretization, let the trial space $\mathcal{U}_h$ be derived from an exact sequence of order $p$. This implies that, at each element $K \in \mathcal{T}$, it is the pullback of $W^p \times \boldsymbol{V}^p|_{\partial\widehat{K}} \cdot \hat{\boldsymbol{n}}_{\widehat{K}}$ from the master square $\widehat{K}$ to the physical element $K$. For the test space $\mathcal{V}_r$, it is sufficient to be from a $p$-enriched sequence of order $p + \Delta p$, locally. Therefore, at each $K$ the test space is homeomorphic to $W^{p+\Delta p}$. Here, however, conformity requirements are not imposed across $\partial K$. That is, $\mathcal{V}_r$ is a space of *discontinuous* piecewise-defined polynomials.

Observe that (5.5) is similar to the standard Bubnov-Galerkin problem; in which, the aim is to find $u \in H_0^1(\Omega) = \mathcal{U}^{\text{BG}}$ such that

$$b^{\text{BG}}(u, v) = \ell^{\text{BG}}(v)\,, \quad \text{for all } v \in H_0^1(\Omega) = \mathcal{V}^{\text{BG}} = \mathcal{U}^{\text{BG}}\,,$$

where

$$\begin{aligned} b^{\text{BG}}(u, v) &= (\operatorname{grad} u, \operatorname{grad} v)_{\Omega}\,, \\ \ell^{\text{BG}}(v) &= (f, v)_{\Omega} - (\operatorname{grad} u^{\text{lift}}, \operatorname{grad} v)_{\Omega}\,. \end{aligned} \tag{5.6}$$

Here, the polynomial space $W^p$ is commonly used for both $\mathcal{U}_h^{\text{BG}} = \mathcal{V}_r^{\text{BG}}$, locally, and continuity is required along the entirety of every element interface.

As far as we are aware, the condition number of the DPG (or DLS) stiffness matrix $\mathsf{A}$ coming from (5.5) has never been derived analytically. We leave that work to another researcher and so do not derive it here, either. Nevertheless, we expect it to grow like $h^{-2}$, as the Bubnov-Galerkin stiffness matrix, $\mathsf{A}_{\mathrm{BG}}$, does [13]. This hypothesis was confirmed with experiments on a square domain $\Omega = (0,1)^2$ starting from a uniform mesh of four square elements and exact solution $u(x,y) = \sin(10\pi x)\sin(10\pi y)$ (see Figure 5.3). Likewise, similar to the least-squares scenario, $\mathrm{cond}(\widetilde{\mathsf{B}}) = \mathrm{cond}(\mathsf{A})^{1/2}$ was confirmed to grow only as $h^{-1}$, and so it eventually became less than the condition number of $\mathsf{A}_{\mathrm{BG}}$, which started out the smallest. Notably, in contrast to least-squares finite element methods, this shows that a first-order system formulation is not required to achieve $\mathcal{O}(h^{-2})$ (or even $\mathcal{O}(h^{-1})$) condition number growth with a DLS method.

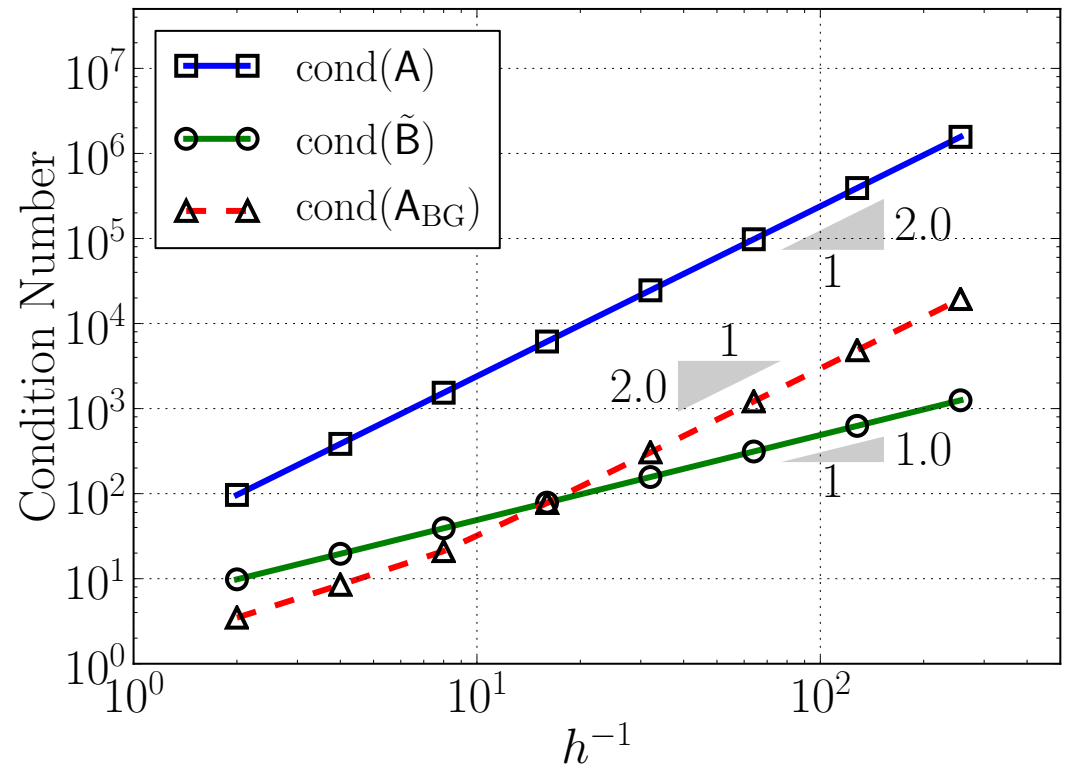

FIGURE 5.3. Comparison with Poisson's equation of the condition number of the Bubnov-Galerkin stiffness matrix $\mathsf{A}_{\mathrm{BG}}$ to the condition number of the DPG stiffness matrices $\mathsf{A}$ and $\widetilde{\mathsf{B}}$ coming from the broken primal formulation, (5.5). Exact sequence polynomial order $p = 2$ for $\mathcal{U}_h$ and $\mathcal{U}_h^{\mathrm{BG}}$, and, in the DPG setting, $\Delta p = 1$. Observe that $\mathrm{cond}(\mathsf{A}_{\mathrm{BG}}) \neq \mathrm{cond}(\mathsf{A})$ and that $\mathrm{cond}(\widetilde{\mathsf{B}}) < \mathrm{cond}(\mathsf{A}_{\mathrm{BG}})$, eventually, for small enough $h$. All reported results are for the statically condensed and diagonally preconditioned matrices.

The optimal solution of the discrete minimum residual problem (2.6) coming from the variational formulation (5.5) would be $\mathbf{u}_{h,\Delta p}^{\mathrm{opt}} = (u_{h,\Delta p}^{\mathrm{opt}}, (\hat{\sigma}_n)_{h,\Delta p}^{\mathrm{opt}})$. If we define the Bubnov-Galerkin problem's solution to be $u^{\mathrm{BG}}$, and the solution of (5.5) to be $(u^{\mathrm{opt}}, \hat{\sigma}_n^{\mathrm{opt}})$ then, it can be shown that $u^{\mathrm{BG}} = u^{\mathrm{opt}}$ [19]. However, for any given $\Delta p$, there is no reason to expect $u_h^{\mathrm{BG}}$ to be equal to $u_{h,\Delta p}^{\mathrm{opt}}$. Indeed, these two solutions did not always agree as is demonstrated in Figure 5.4. Nevertheless, the two different solutions clearly converged to each other, rapidly.

Lastly, as suggested in Section 2.3, it is worth mentioning that, if $f = 0$, the Bubnov-Galerkin solution $u_h^{\mathrm{BG}}$ can also be derived from a minimization principle like (2.11):

$$u_h^{\mathrm{BG}} = \underset{u_h \in \mathcal{U}_h^{\mathrm{BG}}}{\arg\min} \| \operatorname{grad} u_h - \boldsymbol{g} \|_{\boldsymbol{L}^2(\Omega)}^2 , \tag{5.7}$$

where $\boldsymbol{g} = -\operatorname{grad} u^{\mathrm{lift}}$. Likewise, with the DLS theory, it is not difficult to alternatively construct the Bubnov-Galerkin stiffness matrix through a factorization $\mathsf{A}_{\mathrm{BG}} = \widetilde{\mathsf{B}}_{\mathrm{BG}}^* \widetilde{\mathsf{B}}_{\mathrm{BG}}$ [73]. Since this factorization is independent of $f$, this obviously appears appealing because of the reduced condition number of $\widetilde{\mathsf{B}}_{\mathrm{BG}}$. Indeed, had $\mathrm{cond}(\widetilde{\mathsf{B}}_{\mathrm{BG}})$ been computed and added to Figure 5.3 we would certainly verify that $\mathrm{cond}(\widetilde{\mathsf{B}}_{\mathrm{BG}}) = \mathcal{O}(h^{-1})$. Unfortunately, only in the scenario of (5.7) (i.e. when $f = 0$) can orthogonal decomposition methods be used to solve for the solution $u_h^{\mathrm{BG}}$ more accurately and therefore to reduce its sensitivity to round-off error. Furthermore, Bubnov-Galerkin stiffness matrices for most difficult problems will not even be Hermitian/symmetric, so they can not admit *any* such factorization, $\mathsf{A}_{\mathrm{BG}} \neq \mathsf{M}^*\mathsf{M}$ for any matrix $\mathsf{M}$. Therefore, this observation is of extremely limited practical interest.

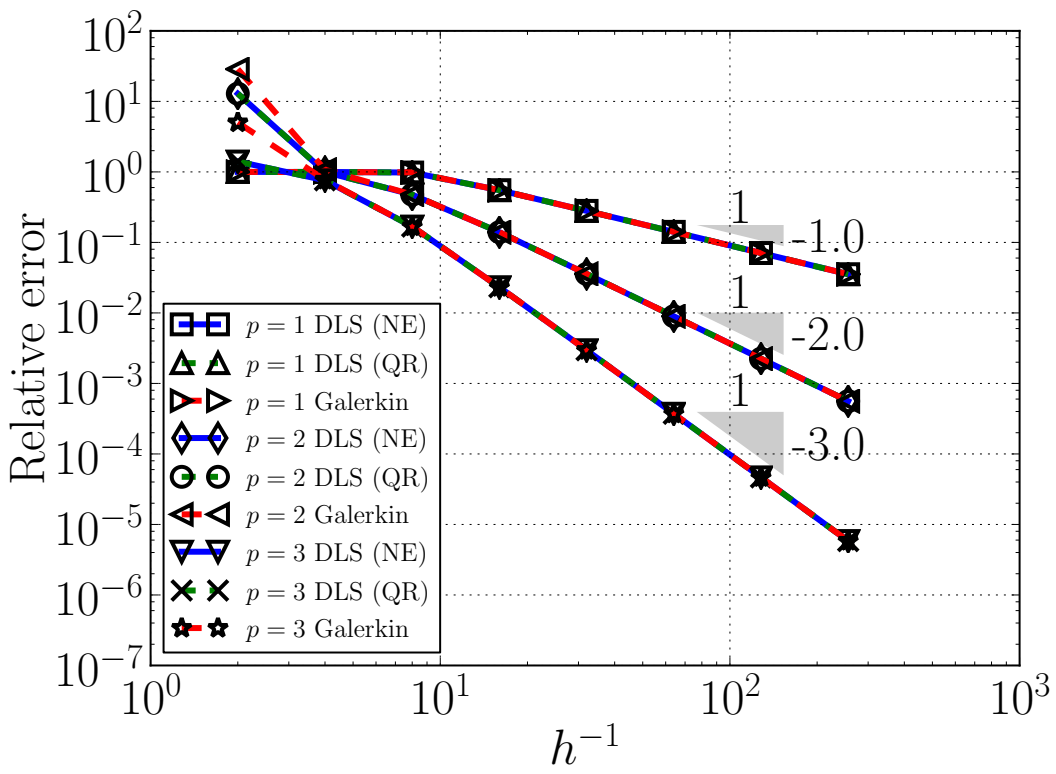

FIGURE 5.4. Relative error in the discrete solutions $u^{\mathrm{opt}}_{h,\Delta p}$ and $u^{\mathrm{BG}}_h$. In the DPG solution, $\Delta p = 1$. Observe that, up to the smallest mesh size considered, the solutions $u^{\mathrm{opt}}_{h,\Delta p}$ obtained by the normal equation (NE) do not differ noticeably from those computed with QR factorization of the least-squares problem (LS).

*Ultraweak DPG methods.* Lastly, consider the broken ultraweak formulation of the Poisson equation in (5.3) (with $\alpha = 0$). This seeks a solution $(u, \boldsymbol{\sigma}, \hat{u}, \hat{\sigma}_n) \in L^2(\Omega) \times \boldsymbol{L}^2(\Omega) \times H^{1/2}_0(\partial\mathcal{T}) \times H^{-1/2}(\partial\mathcal{T}) = \mathfrak{U}$ (again, see [19, 48] for definitions of mesh-trace Sobolev spaces) such that

$$b\big((u, \boldsymbol{\sigma}, \hat{u}, \hat{\sigma}_n), (v, \boldsymbol{\tau})\big) = \ell\big((v, \boldsymbol{\tau})\big), \quad \text{for all } (v, \boldsymbol{\tau}) \in H^1(\mathcal{T}) \times \boldsymbol{H}(\mathrm{div}, \mathcal{T}) = \mathcal{V},$$

where

$$\begin{aligned} b\big((u, \boldsymbol{\sigma}, \hat{u}, \hat{\sigma}_n), (v, \boldsymbol{\tau})\big) &= (\boldsymbol{\sigma}, \mathrm{grad}\, v + \boldsymbol{\tau})_{\mathcal{T}} + (u, \mathrm{div}\,\boldsymbol{\tau})_{\mathcal{T}} - \langle \hat{\sigma}_n, v\rangle_{\partial\mathcal{T}} - \langle \hat{u}, \boldsymbol{\tau}\cdot\hat{\boldsymbol{n}}\rangle_{\partial\mathcal{T}}, \\ \ell\big((v, \boldsymbol{\tau})\big) &= (f, v)_{\mathcal{T}} + \langle \hat{u}^{\mathrm{lift}}, \boldsymbol{\tau}\cdot\hat{\boldsymbol{n}}\rangle_{\partial\mathcal{T}}. \end{aligned} \tag{5.8}$$

Again, $H^1(\mathcal{T})$ and $\boldsymbol{H}(\mathrm{div}, \mathcal{T})$ are broken Sobolev spaces which essentially means that the restriction of each member to any single $K \in \mathcal{T}$ is in $H^1(K)$ and $\boldsymbol{H}(\mathrm{div}, K)$, respectively. We assume that the norm for $\mathcal{V}$ is then $\|\mathfrak{v}\|^2_{\mathcal{V}} = \|v\|^2_{H^1(\mathcal{T})} + \|\boldsymbol{\sigma}\|^2_{\boldsymbol{H}(\mathrm{div},\mathcal{T})}$, although other choices are possible [21]. Lastly, $\hat{u}^{\mathrm{lift}} \in H^{1/2}(\partial\mathcal{T})$ is an extension of $\hat{f} \in H^{1/2}(\partial\Omega)$ to $\partial\mathcal{T}$.

As before, consider a discrete space $\mathfrak{U}_h$ derived from a sequence of order $p$, which would imply that it is (locally) pullbacks of $Y^p \times (Y^p)^2 \times W^p|_{\partial\widehat{K}} \times \boldsymbol{V}^p|_{\partial\widehat{K}}\cdot\hat{\boldsymbol{n}}_{\widehat{K}}$ from the master element $\widehat{K}$ to the physical element $K \in \mathcal{T}$. Here, obtaining a discrete solution $\mathfrak{u}^{\mathrm{opt}}_h$ minimizing the continuous residual (2.4) is not possible and so we have no reference discrete solution to compare different $\mathfrak{u}^{\mathrm{opt}}_{h,r}$ with.

Let the enriched test space $\mathcal{V}_r$ come from a sequence of order $p + \Delta p$, so that locally at each $K$ the space is homeomorphic to $W^{p+\Delta p} \times \boldsymbol{V}^{p+\Delta p}$. In this situation, it was proven in [45] that, provided $\Delta p \geq 2$, the condition number of $\mathsf{A}$ would grow quadratically with $h^{-1}$. However, this was shown only with triangles, where $W^p = \mathcal{P}^p$, $\boldsymbol{V}^p = (\mathcal{P}^{p-1})^2 + x\mathcal{P}^{p-1}$ and $Y^p = \mathcal{P}^{p-1}$, and we have considered quadrilateral elements. Nevertheless, this $\mathcal{O}(h^{-2})$ growth was indeed confirmed with computations as can be observed in Figure 5.5. More importantly, the enriched stiffness matrix $\widetilde{\mathsf{B}}$ was, therefore, verified to have a much improved condition number growth of $\mathcal{O}(h^{-1})$.

5.2. **Failure study.** In some circumstances, finite element stiffness matrices can be so poorly conditioned that the round-off error in solving the discrete equations will compete with, or even surpass, the truncation error coming from the method itself and interpolation spaces being used. In such scenarios, exploiting the overdetermined system of equations with QR is very attractive.

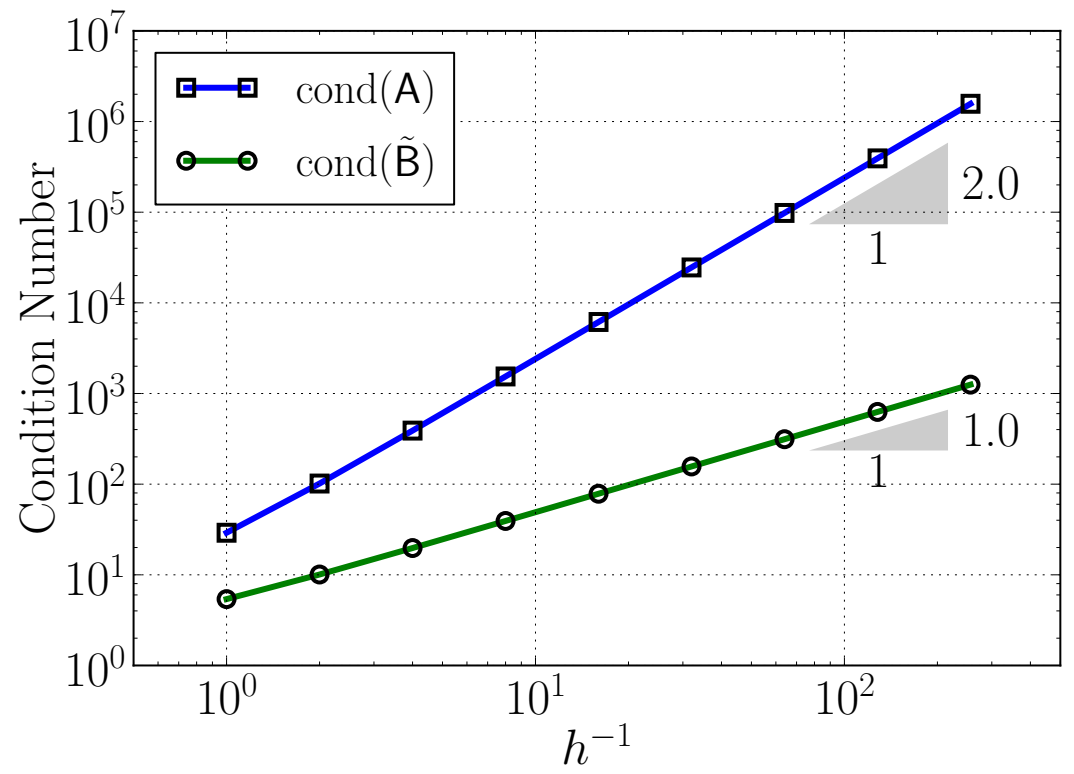

FIGURE 5.5. Condition number growth of the DPG stiffness matrices A and $\widetilde{\mathsf{B}}$ coming from broken ultraweak formulation (5.8). Here, $p = 2$ and $\Delta p = 1$. All reported results are for the statically condensed and diagonally preconditioned matrices.

Due to time and space limitation, we will illustrate this behavior only for ultraweak DPG methods for problems of the form (5.3). We have chosen the ultraweak setting because it is the most actively researched variational setting for DPG methods, at this time. For the normal equation, we solved the system with MUMPS 5.0.1 [52, 1]; and for the overdetermined system, we used qr_mumps 1.2 [15]. Our results are reported in Figure 5.6.

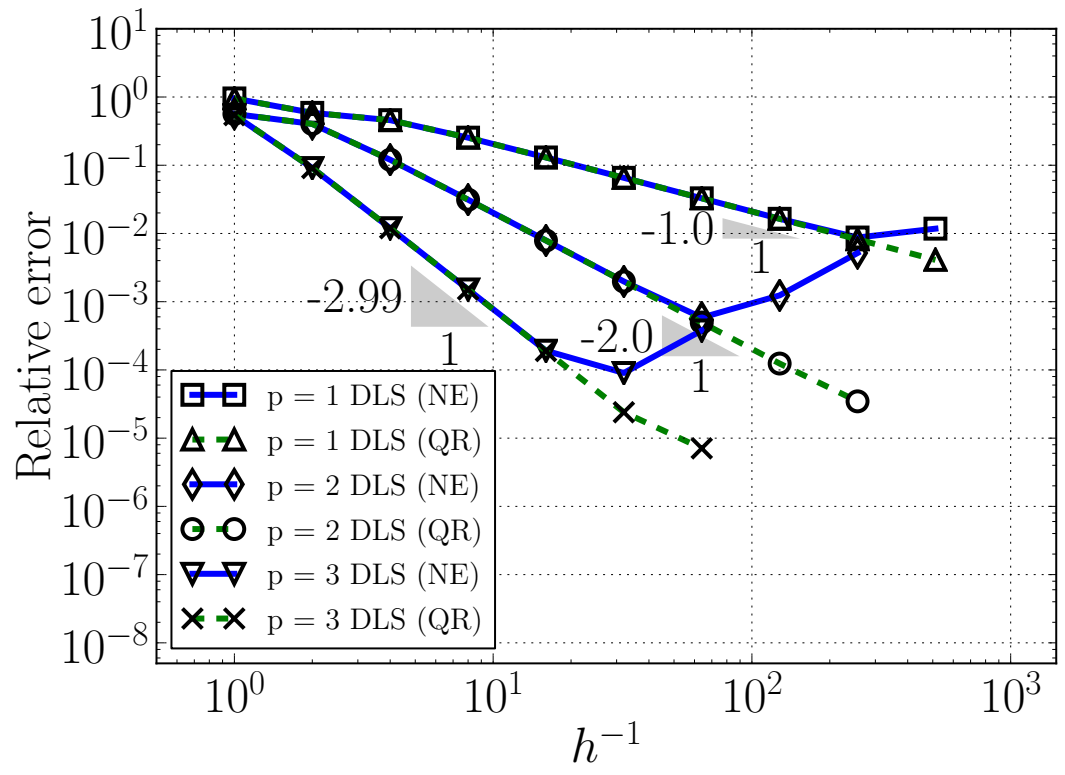

(A) Poisson's equation (single precision).

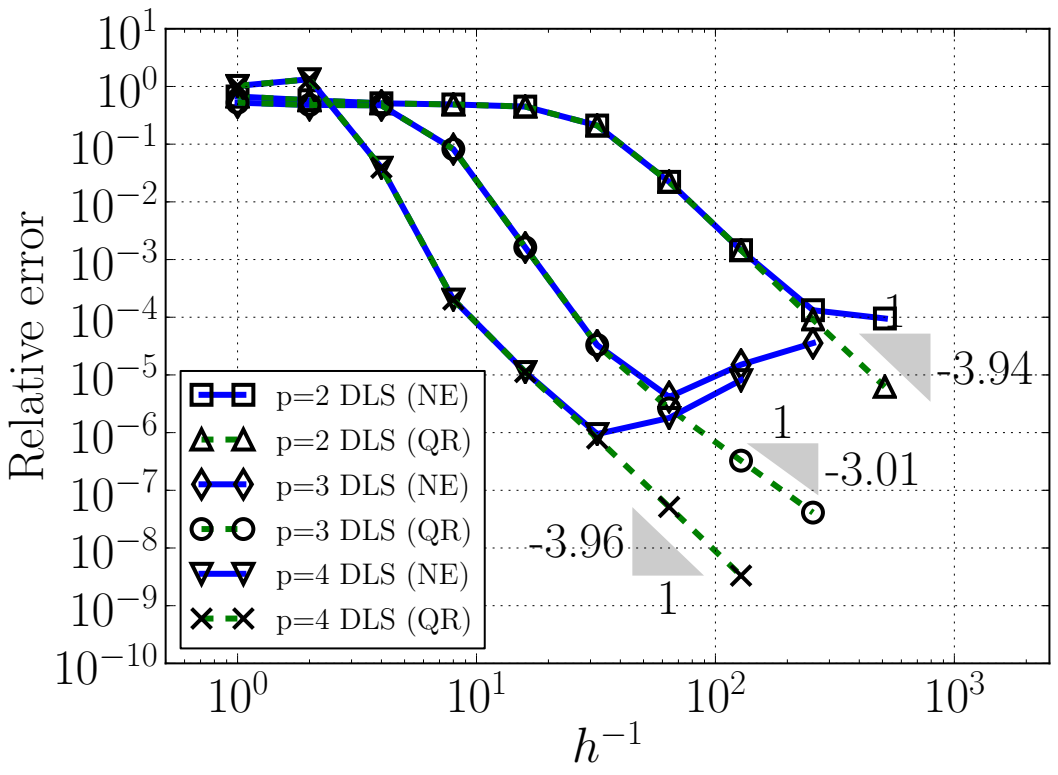

(B) Linear acoustics (double precision).

FIGURE 5.6. Divergence of the discrete solution is observed for various polynomial orders $p$ in two standard ultraweak DPG methods when the normal equation (NE) is constructed, statically condensed, and then solved. Notice that, instead, when QR factorization was used to directly solve the (statically condensed) least-squares problem (LS), the convergence of the discrete solution was maintained, at least for longer. In the $p = 2$ run for (B), had more refinements been performed, we expect that the anticipated rate of convergence (i.e. $2 \neq 3.94$) would have been recovered. $\Delta p = 1$ in all experiments.

First, we performed a *single-precision* floating-point computation with Poisson's equation (5.3) to verify that round-off error would eventually overwhelm truncation error, even in the most well-behaved of problems. Here, we chose $\Omega = (0, 1)^2$ and exact solution $u(x, y) = x^2(1 - x)^2 y^2(1 - y)^2$, and imposed (homogeneous)

Dirichlet boundary conditions around the entire boundary $\partial\Omega$. Uniform $h$-refinements of quadrilateral elements, starting from a single element, were then performed until our computer ran out of memory. In Figure 5.6 (A), we report the loss of convergence to this polynomial exact solution for $p = 1, 2$, and 3 encountered after several mesh refinements when solving the normal equation. Notice that, however, for each polynomial order, qr_mumps applied to the overdetermined system continued to produce the expected rates of convergence after the normal equation approach failed.

In the $p = 2$ case, we can use Figure 5.5 to corroborate this outcome. For example, note that, the two solutions began to visibly diverge at the $7^{\text{th}}$ mesh. Here, the relative error was just below $10^{-3}$. Also, note that, from inspecting Figure 5.5, the corresponding condition number of the stiffness matrix $\mathsf{A}$ was approximately $10^5$. Since machine single precision is approximately $10^{-7}$, the round-off error would have been, at most, approximately $10^{-2}$. It was, therefore, large enough to compete with, or surpass, the truncation error—and this is indeed what we began to see.

Our second example is with the linear acoustics problem in $\Omega = (0, 1)^2$, which we chose to solve with the ultraweak formulation with broken test spaces, *near resonance*, in *double precision*. The equations of linear acoustics are:

$$
\begin{aligned}
i\omega p + \operatorname{div} \mathbf{u} &= f\,, \\
i\omega \mathbf{u} + \operatorname{grad} p &= 0\,.
\end{aligned}
\tag{5.9}
$$

where $p$ is the pressure and $\mathbf{u}$ is the velocity. We imposed a hard boundary for this problem. That is Neumann-type boundary conditions $\mathbf{u}|_{\partial\Omega} \cdot \hat{\boldsymbol{n}} = \hat{g}$ were specified on the entire boundary.

Similar to (5.8), we can derive the corresponding ultraweak formulation. In it, we seek a solution $(p, \mathbf{u}, \hat{p}, \hat{u}_n) \in L^2(\Omega) \times \boldsymbol{L}^2(\Omega) \times H^{1/2}(\partial\mathcal{T}) \times H_0^{-1/2}(\partial\mathcal{T}) = \mathcal{U}$ such that

$$
b\big((p, \mathbf{u}, \hat{p}, \hat{u}_n), (q, \mathbf{v})\big) = \ell\big((q, \mathbf{v})\big)\,, \quad \text{for all } (q, \mathbf{v}) \in H^1(\mathcal{T}) \times \boldsymbol{H}(\operatorname{div}, \mathcal{T}) = \mathcal{V}\,,
$$

where

$$
\begin{aligned}
b\big((p, \mathbf{u}, \hat{p}, \hat{u}_n), (q, \mathbf{v})\big) &= -(p, i\omega q + \operatorname{div} \mathbf{v}) - (\mathbf{u}, i\omega \mathbf{v} + \operatorname{grad} q) + \langle \hat{u}_n, q\rangle_{\partial\mathcal{T}} + \langle \hat{p}, \mathbf{v}\cdot\hat{\boldsymbol{n}}\rangle_{\partial\mathcal{T}}\,, \\
\ell\big((q, \mathbf{v})\big) &= (f, q)_{\mathcal{T}} - \langle \hat{u}_n^{\text{lift}}, q\rangle_{\partial\mathcal{T}}\,.
\end{aligned}
\tag{5.10}
$$

And, where $\hat{u}_n^{\text{lift}} \in H^{-1/2}(\partial\mathcal{T})$ is an extension of $\hat{g} \in H^{-1/2}(\partial\Omega)$ to $\partial\mathcal{T}$. The discretization is also similar to (5.8)—except in $\mathbb{F} = \mathbb{C}$ instead of $\mathbb{R}$—and, again, uses pullbacks of $Y^p \times (Y^p)^2 \times W^p|_{\partial\widehat{K}} \times \boldsymbol{V}^p|_{\partial\widehat{K}}\cdot\hat{\boldsymbol{n}}_{\widehat{K}}$.

We chose to solve the problem for $\omega = 0.5001 \cdot 2\pi$ which is very close to a resonance frequency. Note that, the first eigenvalue of the Laplacian in this setting is $\omega_1 = \pi$. Therefore, we can expect that the stiffness matrix will become very badly-conditioned as the mesh is refined. Using a Gaussian beam for the exact solution (see [61] for equation) and a corresponding discrete lift of $\hat{u}_n^{\text{lift}}$ constructed using projection-based interpolation [27], we clearly received more robust convergence with the QR approach. The results are presented in Figure 5.6 (B).

## 6. Future directions

Many different research pathways for the DLS procedure described here are immediately possible, not the least of which is the exploration of iterative solvers. Another instinctive direction is towards residual problems,

$$
\mathbf{u}_h^{\text{opt}} = \underset{\mathbf{u}_h \in \mathcal{U}_h}{\arg\min} \|\mathcal{F}(\mathbf{u}_h)\|_{\mathcal{V}'}\,,
$$

where $\mathcal{F} : \mathcal{U} \to \mathcal{V}'$ is an appropriate nonlinear operator. In the analogous setting for the statistics community, Gauss-Newton algorithms,

$$
\mathbf{u}_h^{\text{opt}} \mapsto \mathbf{u}_h^{\text{opt}} + \underset{\Delta\mathbf{u}_h \in \mathcal{U}_h}{\arg\min} \|D\mathcal{F}(\mathbf{u}_h)\Delta\mathbf{u}_h + \mathcal{F}(\mathbf{u}_h)\|_{\mathcal{V}'} \qquad \texttt{(repeat)}\,,
$$

are often used to find local minima. And as has been noted in [22], similar techniques are very useful for approximating solutions of nonlinear PDE in the minimum residual fashion. As before, solving the associated normal equation will precede a less satisfactory error bound in each search direction, $\Delta\mathbf{u}_h$, than would occur from handling each least-squares problem directly. This is crucial when the discretization of $D\mathcal{F}(\mathbf{u}_h)$ becomes singular or nearly-singular because the corresponding generalized least-squares problem will always permit a solution, although uniqueness of the solution may be lost.

Moving to more general functional settings than Hilbert spaces—such as featured in Banach space minimum residual problems like those proposed in [55]—a natural extension may exist in applying iteratively reweighted least-squares techniques for a (iteratively updated) DLS discretization. What is more, for problems with parameter uncertainty a total least-squares [53] approach to DLS may have many benefits.

## 7. Conclusion

In this article we have presented a general framework for discrete least-squares finite element methods and illustrated features of this special class of methods. We have demonstrated, both analytically and discretely, how to incorporate essential boundary conditions and equality constraints into this framework and discussed different direct methods to solve the discrete equations. Due to space and time restrictions, our treatment of topics such as preconditioners and iterative solvers is extremely limited. Indeed, these topics are also somewhat open-ended and remain subjects of future research.

We have discussed assembly algorithms for two different formulations of DLS problems which can be solved with separate techniques. In our experiments with one construction, the growth of the condition number of the associated stiffness matrix, after only diagonal preconditioning, was demonstrated to be $\mathcal{O}(h^{-1})$ when solving Poisson's equation. In other experiments which compared sensitivity to round-off error, we demonstrated the associated QR-factorization approach is particularly well-suited for ill-conditioned problems. We have also discussed some techniques which may further improve the conditioning of these matrices during assembly as well as static condensation which can improve the efficiency of the associated solvers.

Although the key results we present are directly towards DPG methods, we have maintained a comprehensive and broad perspective in our presentation which has allowed us to make several connections with other methods in the literature.

**Acknowledgements.** This work was partially supported with grants by NSF (DMS-1418822), AFOSR (FA9550-12-1-0484), and ONR (N00014-15-1-2496).

## Appendix A. Constrained minimization

In this appendix, we demonstrate the necessary and sufficient conditions for existence and uniqueness of the general minimum residual principle over convex sets referred to in Section 2.4. We do not present *a priori* estimates for the rate of converge of the discrete solution. For further inquiring in the DPG setting, see [40].

Let $\mathcal{K}_h \subset \mathcal{U}_h \subset \mathcal{U}$ be a non-empty convex set, let $\|\mathfrak{u}\|_{\mathcal{U}}$ be a semi-norm on $\mathcal{U}$, and assume that the continuous bilinear form $b$ satisfies the continuous inf-sup condition

$$\sup_{\mathfrak{v}\in\mathcal{V}} \frac{b(\mathfrak{u},\mathfrak{v})}{\|\mathfrak{v}\|_{\mathcal{V}}} = \|\mathcal{B}\mathfrak{u}\|_{\mathcal{V}'} \geq \gamma\|\mathfrak{u}\|_{\mathcal{U}}\,, \tag{A.1}$$

for some constant $\gamma > 0$ and all $\mathfrak{u} \in \mathcal{U}$. The corresponding minimum residual problem on $\mathcal{K}_h$ is

$$\min_{\mathfrak{u}_h\in\mathcal{K}_h} \|\mathcal{B}\mathfrak{u}_h - \ell\|_{\mathcal{V}'}\,, \tag{A.2}$$

whose solution satisfies

$$\big(\mathcal{B}\mathfrak{u}_h - \ell, \mathcal{B}(\delta\mathfrak{u}_h - \mathfrak{u}_h)\big)_{\mathcal{V}'} \geq 0 \;\text{ for all } \delta\mathfrak{u}_h \in \mathcal{K}_h\,.$$

Or, equivalently,

$$a(\mathfrak{u}_h, \delta\mathfrak{u}_h - \mathfrak{u}_h) \;\geq\; f(\delta\mathfrak{u}_h - \mathfrak{u}_h) \quad \forall\delta\mathfrak{u}_h \in \mathcal{K}_h\,,$$

where $a : \mathcal{U}\times\mathcal{U} \to \mathbb{F}$ and $f : \mathcal{U} \to \mathbb{F}$ are defined in Section 2.2.

**Theorem A.1.** *The operator $\mathcal{A} = \mathcal{B}'\,\mathcal{R}_{\mathcal{V}}^{-1}\,\mathcal{B}$ is monotone and coercive on $\mathcal{K}_h$ with discrete coercivity constant equal to $\gamma^2$, coming from* (A.1).

*Proof.* Observe that for every $\mathfrak{u}_h, \delta\mathfrak{u}_h \in \mathcal{K}_h$,

$$\begin{aligned}\langle \mathcal{A}\delta\mathfrak{u}_h - \mathcal{A}\mathfrak{u}_h, \delta\mathfrak{u}_h - \mathfrak{u}_h\rangle &= a(\delta\mathfrak{u}_h - \mathfrak{u}_h, \delta\mathfrak{u}_h - \mathfrak{u}_h)\\ &= b(\delta\mathfrak{u}_h - \mathfrak{u}_h, \mathcal{R}_{\mathcal{V}}^{-1}\,\mathcal{B}(\delta\mathfrak{u}_h - \mathfrak{u}_h))\\ &= \langle \mathcal{B}(\delta\mathfrak{u}_h - \mathfrak{u}_h), \mathcal{R}_{\mathcal{V}}^{-1}\,\mathcal{B}(\delta\mathfrak{u}_h - \mathfrak{u}_h)\rangle\\ &= \|\,\mathcal{R}_{\mathcal{V}}^{-1}\,\mathcal{B}(\delta\mathfrak{u}_h - \mathfrak{u}_h)\|_{\mathcal{V}}^2\\ &= \|\mathcal{B}(\delta\mathfrak{u}_h - \mathfrak{u}_h)\|_{\mathcal{V}'}^2\\ &\geq \gamma^2\|\delta\mathfrak{u}_h - \mathfrak{u}_h\|_{\mathcal{U}}^2\,,\end{aligned}$$

from which all results follow. □

**Corollary A.1.** *Solutions to* (A.2) *exist.*

**Corollary A.2.** *If* $\|\cdot\|_{\mathcal{U}}$ *separates points in* $\mathcal{K}_h$ *then* $\mathcal{A}$ *is strictly monotone in* $\mathcal{K}_h$ *and the solution of* (A.2) *is unique.*

The proof of each corollary is classic; see [51].

## APPENDIX B. STATIC CONDENSATION

The purpose of this appendix is to show that the static condensation derived in Section 3.4 coincides with the usual static condensation of finite element matrices via Schur complements. That is, that it coincides with the Schur complement of the normal equation, (NE). With this in mind consider the matrix $\mathsf{A} = \widetilde{\mathsf{B}}^*\widetilde{\mathsf{B}}$, where $\widetilde{\mathsf{B}} = \big[\widetilde{\mathsf{B}}_{\text{bubb.}}\ \widetilde{\mathsf{B}}_{\text{interf.}}\big]$, so that the normal equation is rewritten as

$$\begin{bmatrix}\widetilde{\mathsf{B}}^*_{\text{bubb.}}\widetilde{\mathsf{B}}_{\text{bubb.}} & \widetilde{\mathsf{B}}^*_{\text{bubb.}}\widetilde{\mathsf{B}}_{\text{interf.}}\\ \widetilde{\mathsf{B}}^*_{\text{interf.}}\widetilde{\mathsf{B}}_{\text{bubb.}} & \widetilde{\mathsf{B}}^*_{\text{interf.}}\widetilde{\mathsf{B}}_{\text{interf.}}\end{bmatrix}\begin{bmatrix}\mathbf{u}_{\text{bubb.}}\\ \mathbf{u}_{\text{interf.}}\end{bmatrix} = \mathsf{A}\mathbf{u} = \widetilde{\mathsf{B}}^*\,\widetilde{\mathsf{l}} = \begin{bmatrix}\widetilde{\mathsf{B}}^*_{\text{bubb.}}\widetilde{\mathsf{l}}\\ \widetilde{\mathsf{B}}^*_{\text{interf.}}\widetilde{\mathsf{l}}\end{bmatrix}\,. \tag{B.1}$$

Then writing the Schur complement for $\mathbf{u}_{\text{interf.}}$ yields,

$$\begin{aligned}\widetilde{\mathsf{B}}^*_{\text{interf.}}(\mathrm{I} - \mathsf{P}_{\text{bubb.}})\widetilde{\mathsf{B}}_{\text{interf.}}\mathbf{u}_{\text{interf.}} &= \big(\widetilde{\mathsf{B}}^*_{\text{interf.}}\widetilde{\mathsf{B}}_{\text{interf.}} - \widetilde{\mathsf{B}}^*_{\text{interf.}}\widetilde{\mathsf{B}}_{\text{bubb.}}(\widetilde{\mathsf{B}}^*_{\text{bubb.}}\widetilde{\mathsf{B}}_{\text{bubb.}})^{-1}\widetilde{\mathsf{B}}^*_{\text{bubb.}}\widetilde{\mathsf{B}}_{\text{interf.}}\big)\mathbf{u}_{\text{interf.}}\\ &= \widetilde{\mathsf{B}}^*_{\text{interf.}}\widetilde{\mathsf{l}} - \widetilde{\mathsf{B}}^*_{\text{interf.}}\widetilde{\mathsf{B}}_{\text{bubb.}}(\widetilde{\mathsf{B}}^*_{\text{bubb.}}\widetilde{\mathsf{B}}_{\text{bubb.}})^{-1}\widetilde{\mathsf{B}}^*_{\text{bubb.}}\widetilde{\mathsf{l}} = \widetilde{\mathsf{B}}^*_{\text{interf.}}(\mathrm{I} - \mathsf{P}_{\text{bubb.}})\widetilde{\mathsf{l}}\,,\end{aligned} \tag{B.2}$$

where $\mathsf{P}_{\text{bubb.}} = \widetilde{\mathsf{B}}_{\text{bubb.}}\Big(\widetilde{\mathsf{B}}^*_{\text{bubb.}}\widetilde{\mathsf{B}}_{\text{bubb.}}\Big)^{-1}\widetilde{\mathsf{B}}^*_{\text{bubb.}}$. Then, the only observation that remains is that $\mathsf{P}_{\text{bubb.}} = \mathsf{P}^*_{\text{bubb.}}$ is an orthogonal projection and it can easily be verified that $\mathrm{I} - \mathsf{P}_{\text{bubb.}} = (\mathrm{I} - \mathsf{P}_{\text{bubb.}})^2$, so that

$$\big((\mathrm{I} - \mathsf{P}_{\text{bubb.}})\widetilde{\mathsf{B}}_{\text{interf.}}\big)^*(\mathrm{I} - \mathsf{P}_{\text{bubb.}})\widetilde{\mathsf{B}}_{\text{interf.}}\mathbf{u}_{\text{interf.}} = \big((\mathrm{I} - \mathsf{P}_{\text{bubb.}})\widetilde{\mathsf{B}}_{\text{interf.}}\big)^*(\mathrm{I} - \mathsf{P}_{\text{bubb.}})\widetilde{\mathsf{l}}\,. \tag{B.3}$$

This can easily be seen to be equivalent to an overdetermined system exactly the same as (3.20), and this is the desired result.

The Institute for Computational Engineering and Sciences (ICES), The University of Texas at Austin, 201 E 24th St, Austin, TX 78712, USA